\documentclass[12pt]{amsart}
\usepackage{a4}

\DeclareMathAlphabet{\eusm}{U}{}{}{}  
\SetMathAlphabet\eusm{normal}{U}{eus}{m}{n}
\SetMathAlphabet\eusm{bold}{U}{eus}{b}{n}

\DeclareMathAlphabet{\eufrak}{U}{}{}{}  
\SetMathAlphabet\eufrak{normal}{U}{euf}{m}{n}
\SetMathAlphabet\eufrak{bold}{U}{euf}{b}{n}

\newtheorem{theorem}{Theorem}[section]
\newtheorem{proposition}[theorem]{Proposition}
\newtheorem{lemma}[theorem]{Lemma}
\newtheorem{corollary}[theorem]{Corollary}

\theoremstyle{definition}
\newtheorem{definition}[theorem]{Definition}
\newtheorem{example}[theorem]{Example}

\theoremstyle{remark}
\newtheorem{remark}[theorem]{Remark}
\numberwithin{equation}{section}

\begin{document}
\title{Malliavin Calculus and Skorohod Integration for Quantum Stochastic Processes}
\author{Uwe Franz, R{\'e}mi L{\'e}andre, and Ren\'e Schott}
\address{\underline{U.F.:} Institut f{\"u}r Mathematik und Informatik,
Ernst-Moritz-Arndt-Univer\-sit{\"a}t Greifswald, Jahnstra{\ss}e 15a, D-17487
Greifswald, Germany
\newline \indent
\underline{R.L.:} Institut Elie Cartan, BP 239, 
Universit{\'e} H.~Poincar{\'e}-Nancy I,
F-54506 Vand{\oe}uvre-l\`es-Nancy, France
 \newline \indent
\underline{R.S.:} Institut Elie Cartan and LORIA, BP 239, 
Universit{\'e} H.~Poincar{\'e}-Nancy I,
F-54506 Vand{\oe}uvre-l\`es-Nancy, France
}
\begin{abstract}
A derivation operator and a divergence operator are defined on the algebra of bounded operators on the symmetric Fock space over the complexification of a real Hilbert space $\eufrak{h}$ and it is shown that they satisfy similar properties as the derivation and divergence operator on the Wiener space over $\eufrak{h}$. The derivation operator is then used to give sufficient conditions for the existence of smooth Wigner densities for pairs of operators satisfying the canonical commutation relations. For $\eufrak{h}=L^2(\mathbb{R}_+)$, the divergence operator is shown to coincide with the Hudson-Parthasarathy quantum stochastic integral for adapted integrable processes and  with the non-causal quantum stochastic integrals defined by Lindsay and Belavkin for integrable processes.
\end{abstract}
\maketitle

\markboth{Uwe Franz, Remi L{\'e}andre, and Ren\'e Schott}{Malliavin Calculus and Skorohod Integration for Quantum Stochastic Processes}

\section{Introduction}

Infinite-dimensional analysis has a long history: it began in the sixties (work of Gross \cite{gross67}, Hida, Elworthy, Kr\'ee, $\ldots$), but it is Malliavin \cite{malliavin78} who has applied it to diffusions in order to give a probabilistic proof of H\"ormander's theorem. Malliavin's approach needs a heavy functional analysis apparatus, as the Ornstein-Uhlenbeck operator and the definition of suitable Sobolev spaces, where the diffusions belong. Bismut \cite{bismut81} has given a simpler approach based upon a suitable choice of the Girsanov formula, which gives quasi-invariance formulas. These are differentiated, in order to get integration by parts formulas for the diffusions, which where got by Malliavin in another way.

Our goal is to generalize the hypoellipticity result of Malliavin for non-com\-mutative quantum processes, by using Bismut's method, see also \cite{franz+leandre+schott99}. For that we consider the case of a non-commutative Gaussian process, which is the couple of the position and momentum Brownian motions on Fock space, and we consider the vacuum state. We get an algebraic Girsanov formula, which allows to get integration by parts formulas for the Wigner densities associated to the non-commutative processes, when we differentiate. This allows us to show that the Wigner functional has a density which belongs to all Sobolev spaces over $\mathbb{R}^2$. Let us remark that in general the density is not positive.

If we consider the deterministic elements of the underlying Hilbert space of the Fock space, the derivation of the Girsanov formula leads to a gradient operator satisfying some integration by parts formulas. This shows it is closable as it is in classical infinite-dimensional analysis. But in classical infinite-dimensional analysis, especially in order to study the Malliavin matrix of a functional, we need to be able to take the derivation along a random element of the Cameron-Martin space. In the commutative set-up, this does not pose any problem. Here, we have some difficulty, which leads to the definition of a right-sided and a left-sided gradient, which can be combined to a two-sided gradient.

We can define a divergence operator as a kind of adjoint of the two-sided gradient for cylindrical (non-commutative) vector fields, but since the vacuum state does not define a Hilbert space, it is more difficult to extend it to general (non-commutative) vector fields.

We show that the non-commutative differential calculus contains in some sense the commutative differential calculus.

In the white noise case, i.e.\ if the underlying Hilbert space is the $L^2$-space of some measure space, the classical divergence operator defines an anticipating stochastic integral, known as the Hitsuda-Skorohod integral. We compute the matrix elements between exponential vectors for our divergence operator and use them to show that the divergence operator coincides with the non-causal creation and annihilation integrals defined by Belavkin \cite{belavkin91a,belavkin91b} and Lindsay \cite{lindsay93} for integrable processes, and therefore with the Hudson-Parthasarathy \cite{hudson+parthasarathy84} integral for adapted processes.

\section{Analysis on Wiener space}\label{comm Wiener}

Let us first briefly recall a few definitions and facts from analysis on Wiener space, for more details see, e.g., \cite{janson97,malliavin97,nualart95,nualart98,ustunel95}. Let $\eufrak{h}$ be a real separable Hilbert space. Then there exists a probability space $(\Omega,\mathcal{F},\mathbb{P})$ and a linear map $W:\eufrak{h}\to L^2(\Omega)$ such that the $W(h)$ are centered Gaussian random variables with covariances given by
\[
\mathbb{E}\big(W(h)W(k)\big) = \langle h,k\rangle, \qquad \text{ for all } h,k\in\eufrak{h}.
\]
Set $\mathcal{H}_1=W(\eufrak{h})$, this is a closed Gaussian subspace of $L^2(\Omega)$ and $W:\eufrak{h}\to\mathcal{H}_1\subseteq L^2(\Omega)$ is an isometry. We will assume that the $\sigma$-algebra $\mathcal{F}$ is generated by the elements of $\mathcal{H}_1$. We introduce the algebra of bounded smooth functionals
\[
\mathcal{S}=\{F=f\big(W(h_1),\ldots,W(h_n)\big)| n\in\mathbb{N}, f\in C^\infty_b(\mathbb{R}^n), h_1,\ldots,h_n\in\eufrak{h}\},
\]
and define the derivation operator $\tilde{D}:\mathcal{S}\to L^2(\Omega)\otimes \eufrak{h}\cong L^2(\Omega;\eufrak{h})$ by
\[
\tilde{D}F=\sum_{i=1}^n \frac{\partial f}{\partial x_i}\big(W(h_1),\ldots,W(h_n)\big)\otimes h_i
\]
for $F=f\big(W(h_1),\ldots,W(h_n)\big)\in\mathcal{S}$. Then one can verify the following properties of $\tilde{D}$.
\begin{enumerate}
\item
$\tilde{D}$ is a derivation (w.r.t.\ the natural $L^\infty(\Omega)$-bimodule structure of $L^2(\Omega;\eufrak{h})$), i.e.\
\[
\tilde{D}(FG)= F(\tilde{D}G)+(\tilde{D}G)F, \qquad \text{ for all } F,G\in\mathcal{S}.
\]
\item
The scalar product $\langle h, \tilde{D}F\rangle $ coincides with the Fr\'echet derivative
\[
\tilde{D}_hF = \left.\frac{{\rm d}}{{\rm d}\varepsilon}\right|_{\varepsilon=0} f\big(W(h_1)+\varepsilon\langle h,h_1\rangle,\ldots,W(h_n)+\varepsilon\langle h,h_n\rangle\big)
\]
for all $F=f\big(W(h_1),\ldots,W(h_n)\big)\in\mathcal{S}$ and all $h\in\eufrak{h}$.
\item
We have the following integration by parts formulas,
\begin{eqnarray}
\mathbb{E}\big(FW(h)\big) &=& \mathbb{E}\big(\langle h, \tilde{D}F\rangle\big) \label{int by parts1}\\
\mathbb{E}\big(FGW(h)\big) &=& \mathbb{E}\big(\langle h, \tilde{D}F\rangle G+F \langle h, \tilde{D}G\rangle \big) \label{int by parts2}
\end{eqnarray}
for all $F,G\in\mathcal{S}$, $h\in\eufrak{h}$.
\item
The derivation operator $\tilde{D}$ is a closable operator from $L^p(\Omega)$ to $L^p(\Omega;\eufrak{h})$ for $1\le p\le \infty$. We will denote its closure again by $\tilde{D}$.
\end{enumerate}
We can also define the gradient $\tilde{D}_u F= \langle u,DF\rangle$ w.r.t.\ $\eufrak{h}$-valued random variables $u\in L^2(\Omega;\eufrak{h})$, this is $L^\infty(\Omega)$-linear in the first argument and a derivation in the second, i.e.\
\begin{eqnarray*}
\tilde{D}_{Fu}G &=& F \tilde{D}_uG, \\
\tilde{D}_u (FG) &=& F (\tilde{D}_uG)+ (\tilde{D}_uF)G.
\end{eqnarray*}

$L^2(\Omega)$ and $L^2(\Omega;\eufrak{h})$ are Hilbert spaces (with the obvious inner products), therefore the closability of $\tilde{D}$ implies that it has an adjoint. We will call the adjoint of $\tilde{D}:L^2(\Omega)\to L^2(\Omega;\eufrak{h})$ the divergence operator and denote it by $\tilde\delta:L^2(\Omega;\eufrak{h})\to L^2(\Omega)$. Denote by \[
\mathcal{S}_\eufrak{h}=\left\{u=\sum_{j=1}^n F_j \otimes h_j \Big| n\in\mathbb{N},F_1,\ldots,F_n\in \mathcal{S}, h_1,\ldots,h_n\in\eufrak{h}\right\} 
\]
the smooth elementary $\eufrak{h}$-valued random variables, then $\tilde\delta(u)$ is given by
\[
\tilde\delta(u)=\sum_{j+1}^n F_j W(h_j) - \sum_{j+1}^n \langle h_j, \tilde{D}F_j\rangle
\]
for $u=\sum_{j=1}^n F_j\otimes h_j\in\mathcal{S}_\eufrak{h}$. If we take, e.g., $\eufrak{h}=L^2(\mathbb{R}_+)$, then $B_t=W(\mathbf{1}_{[0,t]})$ is a standard Brownian motion, and the $\eufrak{h}$-valued random variables can also be interpreted as stochastic processes indexed by $\mathbb{R}_+$. It can be shown that $\tilde\delta(u)$ coincides with the It\^o integral $\int_{\mathbb{R}_+} u_t{\rm d}W_t$ for adapted integrable processes. In this case the divergence operator is also called the Hitsuda-Skorohod integral.

The derivation operator and the divergence operator satisfy the following relations
\begin{eqnarray}
\tilde{D}_h\big(\delta(u)) &=& \langle h, u\rangle + \tilde\delta(\tilde{D}_hu), \label{comm rel1}\\
\mathbb{E}\big( \tilde\delta(u)\tilde\delta(v)\big) &=& \mathbb{E}\big(\langle u,v\rangle \big) + \mathbb{E}\big( {\rm Tr}(\tilde{D}u\circ \tilde{D}v)\big), \label{skorohod cov1}\\
\tilde\delta(Fu) &=& F \tilde\delta(u) - \langle u, \tilde{D}F\rangle, \label{DFu1}
\end{eqnarray}
for $h\in\eufrak{h}$, $u,v\in\mathcal{S}_\eufrak{h}$, $F\in\mathcal{S}$.
Here $\tilde{D}$ is extended in the obvious way to $\eufrak{h}$-valued random variables, i.e.\ as $\tilde{D}\otimes{\rm id}_\eufrak{h}$. Thus $\tilde{D}u$ is an $\eufrak{h}\otimes\eufrak{h}$-valued random variable and can also be interpreted as a random variable whose values are (Hilbert-Schmidt) operators on $\eufrak{h}$. If $\{e_j;j\in\mathbb{N}\}$ is a complete orthonormal system on $\eufrak{h}$, then ${\rm Tr}(\tilde{D}u\circ \tilde{D}v)$ can be computed as
${\rm Tr}(\tilde{D}u\circ \tilde{D}v) = \sum_{i,j=1}^\infty \tilde{D}_{e_i}\langle u, e_j\rangle \tilde{D}_{e_j}\langle v,e_i\rangle$.

\section{The non-commutative Wiener space}\label{non comm Wiener}

Let again $\eufrak{h}$ be a real separable Hilbert space and let $\eufrak{h}_\mathbb{C}$ be its complexification. Then we can define a conjugation $\overline{\phantom{a}}:\eufrak{h}_\mathbb{C}\to \eufrak{h}_\mathbb{C}$ by $\overline{h_1+ih_2}=h_1-ih_2$ for $h_1,h_2\in\eufrak{h}_\mathbb{C}$. This conjugation satisfies $\left\langle \overline{h},\overline{k}\right\rangle = \overline{\langle h,k\rangle} = \langle k,h\rangle$ for all $h,k\in\eufrak{h}_\mathbb{C}$. The elements of $\eufrak{h}$ are characterized by the property $\overline{h}=h$, we will call them real. 

Let $\eufrak{H}=\Gamma_s(\eufrak{h}_\mathbb{C})$ be the symmetric Fock space over $\eufrak{h}_\mathbb{C}$, i.e.\ $\eufrak{H}= \bigoplus_{n\in\mathbb{N}} \eufrak{h}_\mathbb{C}^{\odot n}$, where `$\odot$' denotes the symmetric tensor product, and denote the vacuum vector $1+0+\cdots$ by $\Omega$. It is well-known that the symmetric Fock space is isomorphic to the complexification of the Wiener space $L^2(\Omega)$ associated to $\eufrak{h}$ in Section \ref{comm Wiener}. We will develop a calculus on the non-commutative probability space $(\mathcal{B}(\eufrak{H}),\mathbb{E})$, where $\mathbb{E}$ denotes the state defined by $\mathbb{E}(X)=\langle \Omega,X\Omega\rangle$ for $X\in \mathcal{B}(\eufrak{H})$. To emphasize the analogy with the analysis on Wiener space we call $(\mathcal{B}(\eufrak{H}),\mathbb{E})$ the non-commutative Wiener space over $\eufrak{h}$.

The exponential vectors $\{\mathcal{E}(k)= \sum_{n=0}^\infty \frac{k^{\otimes n}}{\sqrt{n!}}; k\in\eufrak{h}_\mathbb{C}\}$ are total in $\eufrak{H}$, their scalar product is given by
\[
\langle\mathcal{E}(k_1),\mathcal{E}(k_2)\rangle = e^{\langle k_1,k_2\rangle}.
\]

We can define the operators $a(h),a^+(h),Q(h),P(h)$ (annihilation, creation, position, momentum) and $U(h_1,h_2)$ with $h,h_1,h_2\in\eufrak{h}_\mathbb{C}$ on $\eufrak{H}$, see, e.g., \cite{biane93,meyer95,parthasarathy92}. The creation and annihilation operators $a^+(h)$ and $a(h)$ are closed, unbounded, mutually adjoint operators. The position and momentum operators
\[
Q(h)= \big(a(\overline{h})+a^+(h)\big), \quad \text{ and }\quad
P(h)= i\big(a(\overline{h})-a^+(h)\big)
\]
are self-adjoint, if $h$ is real. 

The commutation relations of creation, annihilation, position, and momentum are
\[
\begin{array}{rclcrcl}
[a(h),a^+(k)]&=& \langle h,k\rangle, &&
[a(h),a(k)]&=&[a^+(h),a^+(k)]=0, \\ {}
[Q(h),Q(k)]&=& [P(h),P(k)] = 0, &&
[P(h),Q(k)]&=& 2 i \langle \overline{h},k\rangle.
\end{array}
\]
The Weyl operators $U(h_1,h_2)$ can be defined
by $U(h_1,h_2)= \exp\big(iP(h_1)+iQ(h_2)\big)=\exp i \big(a(\overline{h_2}-i\overline{h_1})+ a^+(h_2-ih_1)\big)$, they satisfy
\[
U(h_1,h_2) U(k_1,k_2) = \exp i\big( \langle \overline{h}_2,k_1\rangle - \langle \overline{h}_1,k_2\rangle\big) U(h_1+h_2,k_1+k_2)
\]
Furthermore we have $U(h_1,h_2)^*=U(-\overline{h}_1,-\overline{h}_2)$ and $U(h_1,h_2)^{-1}=U(-h_1,-h_2)$. We see that $U(h_1,h_2)$ is unitary, if $h_1$ and $h_2$ are real. These operators act on the vacuum $\Omega=\mathcal{E}(0)$ as
\[
U(h_1,h_2)\Omega = \exp\left( -\frac{\langle\overline{h}_1,h_1\rangle +\langle\overline{h}_2,h_2\rangle}{2} \right)\mathcal{E}\left(h_1+ih_2\right)
\]
and on general exponential vectors $\mathcal{E}\left(f\right)= \sum_{n=0}^\infty \frac{f^{\otimes n}}{\sqrt{n!}}$ as
\[
U(h_1,h_2)\mathcal{E}\left(f\right) = \exp\left(-\langle \overline{f},h_1+ih_2\rangle -\frac{\langle\overline{h}_1,h_1\rangle +\langle\overline{h}_2,h_2\rangle}{2}\right)\mathcal{E}\left(f+h_1+ih_2\right).
\]
The operators $a(h),a^+(h),Q(h),P(h)$ and $U(h_1,h_2)$ are unbounded, but their domains contain the exponential vectors. We will want to compose them with bounded operators on $\eufrak{H}$, to do so we adopt the following convention. Let
\begin{eqnarray*}
\mathcal{L}\big(\mathcal{E}(\eufrak{h}_\mathbb{C}),\eufrak{H}\big)&=&\Big\{B\in{\rm Lin}\,\big({\rm span}(\mathcal{E}(\eufrak{h}_\mathbb{C})),\eufrak{H}\big)\Big| \exists B^* \in{\rm Lin}\,\big({\rm span}(\mathcal{E}(\eufrak{h}_\mathbb{C})),\eufrak{H}\big)\\
 &&\quad\text{ s.t. } \big\langle\mathcal{E}(f),B\mathcal{E}(g)\big\rangle = \big\langle B^*\mathcal{E}(f),\mathcal{E}(g)\big\rangle \text{ for all }f,g\in\eufrak{h}_\mathbb{C}\Big\},
\end{eqnarray*}
i.e.\ the space of linear operators that are defined on the exponential vectors and that have an ``adjoint'' that is also defined on the exponential vectors. Obviously $a(h),a^+(h),Q(h),P(h), U(h_1,h_2)\in \mathcal{L}\big(\mathcal{E}(\eufrak{h}_\mathbb{C}),\eufrak{H}\big)$. We will say that an expression of the form
$\sum_{j=1}^n X_j B_j Y_j$
with $X_1,\ldots,X_n,Y_1,\ldots,Y_n\in\mathcal{L}\big(\mathcal{E}(\eufrak{h}_\mathbb{C}),\eufrak{H})$ and $B_1,\ldots,B_n\in \mathcal{B}(\eufrak{H})$ defines a bounded operator on $\eufrak{H}$, if there exists a bounded operator $M\in \mathcal{B}(\eufrak{H})$ such that
\[
\big\langle\mathcal{E}(f),M\mathcal{E}(g)\big\rangle = \sum_{j=1}^n\big\langle X_j^*\mathcal{E}(f),B_jY_j\mathcal{E}(g)\big\rangle
\]
holds for all $f,g\in\eufrak{h}_\mathbb{C}$. If it exists, this operator is unique, because the exponential vectors are total in $\eufrak{H}$. We will then write
\[
M=\sum_{j=1}^n X_j B_j Y_j.
\]

\section{Weyl calculus}

\begin{definition}\label{def-weyl}
Let $h=(h_1,h_2)\in \eufrak{h}\otimes \mathbb{R}^2$. We set
\begin{eqnarray}
{\rm Dom}\, O_h &=& \Big\{\varphi:\mathbb{R}^2\to\mathbb{C}\Big|\exists M\in\mathcal{B}(\eufrak{H}),\forall k_1,k_2\in\eufrak{h}_\mathbb{C}:\langle\mathcal{E}(k_1),M\mathcal{E}(k_2)\rangle = \nonumber \\
&&\qquad\qquad\frac{1}{2\pi}\int\langle\mathcal{E}(k_1),U(uh_1,vh_2)\mathcal{E}(k_2)\rangle \mathcal{F}^{-1}\varphi(u,v){\rm d}u{\rm d}v\Big\} \label{defin-O_h}
\end{eqnarray}
and for $\varphi\in{\rm Dom}\,O_h$ we define $O_h(\varphi)$ to be the bounded operator $M$ appearing in Equation (\ref{defin-O_h}), it is uniquely determined due to the totality of $\{\mathcal{E}(k): k\in\eufrak{h}_\mathbb{C}\}$.
\end{definition}
We take the Fourier transform $\mathcal{F}$ as
\[
\mathcal{F}\varphi(u,v)= \frac{1}{2\pi} \int_{\mathbb{R}^2} \varphi(x,y) \exp\big(i(ux+vy)\big) {\rm d}x{\rm d}y.
\]
Its inverse is simply
\[
\mathcal{F}^{-1}\varphi(x,y)= \frac{1}{2\pi} \int_{\mathbb{R}^2} \varphi(u,v) \exp\big(-i(ux+vy)\big) {\rm d}u{\rm d}v.
\]

\begin{remark}
If $\varphi$ is a Schwartz function on $\mathbb{R}^2$, then one can check that $O_h(\varphi) = \frac{1}{2\pi}\int_{\mathbb{R}^{2}} \mathcal{F}^{-1}\varphi(u,v) \exp\big( iuP(h_1)+iv Q(h_2)\big) {\rm d}u{\rm d}v$ defines a bounded operator. It is known that the map from $\mathcal{S}(\mathbb{R}^2)$ to $B(\eufrak{H})$ defined in this way extends to a continuous map from $L^p(\mathbb{R})$ to $B(\eufrak{H})$ for all $p\in[1,2]$, but that for $p>2$ there exist functions in $L^p(\mathbb{R}^2)$ for which we can not define a bounded operator in this way, see, e.g., \cite{wong98} and the references cited therein. But it can be extended to exponential functions, since $\frac{1}{2\pi}\mathcal{F}^{-1}\exp i(x_0u+y_0v)= \delta_{(x_0,y_0)}$ and thus
\[
O_h\big(\exp i(x_0u+y_0v)\big)=U(x_0h_1,y_0h_2).
\]
\end{remark}

\begin{lemma}\label{lem-bounded}
Let $1\le p\le 2$ and $h\in\eufrak{h}\otimes\mathbb{R}^2$ such that $\langle h_1,h_2\rangle\not=0$. Then we have $L^p(\mathbb{R}^2)\subseteq{\rm Dom}\,O_h$ and there exists a constant $C_{h,p}$ such that
\[
||O_h(\varphi)||\le C_{h,p} ||\varphi||_p
\]
for all $\varphi\in L^p(\mathbb{R}^2)$.
\end{lemma}
\begin{proof}
This follows immediately from \cite[Theorem 11.1]{wong98}, where it is stated for the irreducible unitary representation with parameter $\hbar=1$ of the Heisenberg-Weyl group.
\end{proof}

As `joint density' of the pair $\big(P(h_1),Q(h_2)\big)$ we will use its Wigner distribution.
\begin{definition}
Let $\Phi$ be a state on $B(\eufrak{H})$. We will call ${\rm d}W_{h,\Phi}$ the Wigner distribution of $\big(P(h_1),Q(h_2)\big)$ in the state $\Phi$, if
\[
\int \varphi {\rm d}W_{h,\Phi} = \Phi\big(O_h(\varphi)\big)
\]
is satisfied for all Schwartz functions $\varphi$.
\end{definition}

In general, ${\rm d}W_{h,\Phi}$ is not positive, but only a signed measure, since $O_h$ does not map positive functions to positive operators. But we can show that it has a density.
\begin{proposition}
Let $h=(h_1,h_2)\in \eufrak{h}\otimes \mathbb{R}^2$ such that $\langle h_1,h_2\rangle\not=0$ and let $\Phi$ be a state on $B(\eufrak{H})$. Then there exists a function $w_{h,\Phi}\in \bigcap_{2\le p\le\infty} L^p(\mathbb{R}^2)$ such that ${\rm d}W_{h,\Phi}= w_{h,\Phi}{\rm d}x{\rm d}y$.
\end{proposition}
\begin{proof}
It is sufficient to observe that Lemma \ref{lem-bounded} implies that the map $\varphi\mapsto \Phi\big(O_h(\varphi)\big)$ defines a continuous linear functional on $L^p(\mathbb{R}^2)$ for $1\le p\le 2$.
\end{proof}

The following proposition will play the role of the Girsanov transformation in classical Malliavin calculus. If we conjugate $O_h(\varphi)$ with $U(-k_2/2,k_1/2)$ for $k\in \eufrak{h}\otimes\mathbb{R}^2$, then this amounts to a translation of the argument of $\varphi$ by $(\langle k_1,h_1\rangle,\langle k_2,h_2\rangle)$.

\begin{proposition}\label{prop girsanov}
Let $h,k\in \eufrak{h}\otimes\mathbb{R}^2$ and $\varphi\in {\rm Dom}\,O_h$. Then we have
\[
U(-k_2/2,k_1/2)O_h(\varphi)U(-k_2/2,k_1/2)^*= O_h\big(T_{(\langle k_1,h_1\rangle,\langle k_2,h_2\rangle)}\varphi\big)
\]
where $T_{(x_0,y_0)}\varphi(x,y)=\varphi(x+x_0,y+y_0)$.
\end{proposition}
\begin{proof}
For $(u,v)\in\mathbb{R}^2$, we have
\begin{eqnarray*}
&&U(-k_2/2,k_1/2)\exp\big(i(uP(h_1)+vQ(h_2))\big)U(-k_2/2,k_1/2)^* \\
&=& U(-k_2/2,k_1/2)U(uh_1,vh_2)U(-k_2/2,k_1/2)^* \\
&=&\exp -i\big(u\langle k_1, h_1\rangle+v\langle k_2, h_2\rangle\big) U(uh_1,vh_2)
\end{eqnarray*}
and therefore
\begin{eqnarray*}
&& U(-k_2/2,k_1/2)O_h(\varphi)U(-k_2/2,k_1/2)^* \\
&=&\int_{\mathbb{R}^{2}} \mathcal{F}^{-1}\varphi(u,v)\exp\left(-i\big(u\langle k_1, h_1\rangle+v\langle k_2, h_2\rangle \big)\right) \exp i\big(uP(h_1)+vQ(h_2)\big) {\rm d}u{\rm d}v \\
&=&\int_{\mathbb{C}^{d}} \mathcal{F}^{-1}T_{(\langle k_1,h_1\rangle,\langle k_2,h_2\rangle)}  \varphi(u,v) \exp i\big(uP(h_1)+vQ(h_2)\big) {\rm d}u{\rm d}v \\
&=& O_h\left( T_{(\langle k_1,h_1\rangle,\langle k_2,h_2\rangle)}\varphi\right).
\end{eqnarray*}
\end{proof}

{}From this formula we can derive a kind of integration by parts formula that can be used to get the estimates that show the differentiability of the Wigner densities.
\begin{proposition}\label{prop-int by parts}
Let $h\in \eufrak{h}\otimes \mathbb{R}^2$, $k\in \eufrak{h}_\mathbb{C}\otimes \mathbb{C}^2$, and $\varphi$ such that $\varphi,\frac{\partial\varphi}{\partial x},\frac{\partial\varphi}{\partial y}\in {\rm Dom}\,O_h$. Then $[Q(\overline{k}_1)-P(\overline{k}_2),O_h(\varphi)]$ defines a bounded operator on $\eufrak{H}$ and we have
\[
\frac{i}{2}[Q(\overline{k}_1)-P(\overline{k}_2),O_h(\varphi)] = O_h\left(\langle k_1, h_1\rangle\frac{\partial\varphi}{\partial x}+ \langle k_2, h_2\rangle\frac{\partial\varphi}{\partial y}\right)
\]
\end{proposition}
\begin{proof}
For real $k$ this is the infinitesimal version of the previous proposition, just differentiate
\[
U(\varepsilon k_2/2,\varepsilon k_1/2)O_h(\varphi)U(\varepsilon k_2/2,\varepsilon k_1/2)^*=O_h\big(T_{(\varepsilon \langle k_1,h_1\rangle,\varepsilon\langle k_2,h_2\rangle)}\varphi\big)
\]
with respect to $\varepsilon$ and set $\varepsilon =0$. For complex $k$ it follows by linearity.
\end{proof}

Like the integration by parts formula in classical Malliavin calculus, this formula follows from a Girsanov transformation. Furthermore, it can also be used to derive sufficient conditions for the existence of smooth densities.

\begin{proposition}\label{prop-regular1}
Let $\kappa\in\mathbb{N}$, $h\in\eufrak{h}\otimes\mathbb{R}^2$ with $\langle h_1,h_2\rangle\not=0$, and $\Phi$ a vector state, i.e.\ there exists a unit vector $\omega\in\eufrak{H}$ such that $\Phi(X)=\langle\omega,X\omega\rangle$ for all $X\in B(\eufrak{H})$. If there exists a $k\in\eufrak{h}_\mathbb{C}\otimes\mathbb{C}^2$ such that
\[
\omega\in\bigcap_{\kappa_1+\kappa_2\le \kappa}{\rm Dom}\,Q(k_1)^{\kappa_1}P(k_2)^{\kappa_2}\cap \bigcap_{\kappa_1+\kappa_2\le \kappa}{\rm Dom}\,Q(\overline{k}_1)^{\kappa_1}P(\overline{k}_2)^{\kappa_2}
\]
and
\[
\langle h_1,k_1\rangle \not=0\qquad\text{ and }\qquad \langle h_2,k_2\rangle \not=0,
\]
then $w_{h,\Phi}\in \bigcap_{2\le p\le \infty}H^{p,\kappa}(\mathbb{R}^2)$, i.e.\ the Wigner density $w_{h,\Phi}$ lies in the Sobolev spaces of order $\kappa$ for all $2\le p\le \infty$. 
\end{proposition}
\begin{proof}
We will show the result for $\kappa=1$, the general case can be shown similarly (see also the proof of Theorem \ref{theo regular2}). Let $\varphi$ be a Schwartz function. Let $p\in[1,2]$. Then we have
\begin{eqnarray*}
\left|\int\frac{\partial\varphi}{\partial x}{\rm d}W_{h,\Phi}\right| &=& \left|\left\langle \omega, O_h\left(\frac{\partial\varphi}{\partial x}\right)\omega\right\rangle\right| \\
&=& \left|\left\langle \omega,\frac{i}{2|\langle k_1,h_1\rangle|} \big[Q(\overline{k}_1), O_h(\varphi)\big]\omega\right\rangle\right| \\
&\le&
\frac{C_{h,p} \big(||Q(k_1)\omega||+||Q(\overline{k}_1)\omega||\big)}{2|\langle k_1,h_1\rangle|} ||\varphi||_p.
\end{eqnarray*}
Similarly, we get
\[
\left|\int\frac{\partial\varphi}{\partial y}{\rm d}W_{h,\Phi}\right| \le
\frac{C_{h,p} \big(||P(k_2)\omega||+||P(\overline{k}_2)\omega||\big)}{2|\langle k_2,h_2\rangle|} ||\varphi||_p,
\]
and together these two inequalities imply $w_{h,\Phi}\in H^{p',1}(\mathbb{R}^2)$ for $p'=\frac{p}{p-1}$.
\end{proof}
We will give a more general result of this type in Theorem \ref{theo regular2}.

\section{The derivation operator}

In this section we define a derivation operator on our non-commutative probability space and show that it satisfies similar properties as the derivation operator on Wiener space.

We want to interpret the expression in the integration by parts formula in Proposition \ref{prop-int by parts} as a directional or Fr\'echet derivative.
\begin{definition}\label{def derivation1}
Let $k\in \eufrak{h}_\mathbb{C}\otimes\mathbb{C}^2$. We set
\begin{eqnarray*}
{\rm Dom}\,D_k &=&\left\{B\in \mathcal{B}(\eufrak{H})\Big|\,\frac{i}{2}[Q(k_1)-P(k_2),B] \text{ defines a bounded operator on }\eufrak{H}\right\}
\end{eqnarray*}
and for $B\in {\rm Dom}\,D_k$, we set $D_k B=\frac{i}{2}[Q(k_1)-P(k_2),B]$.
\end{definition}

Note that $B\in{\rm Dom}\,D_k$ for some $k\in\eufrak{h}_\mathbb{C}\otimes\mathbb{C}^2$ implies $B^*\in{\rm Dom}\,D_{\overline{k}}$ and
\[
D_{\overline{k}}B^* = (D_kB)^*.
\]

\begin{example}
Let $k\in \eufrak{h}_\mathbb{C}\otimes\mathbb{C}^2$ and let $\psi\in{\rm Dom}\,P(k_2)\cap{\rm Dom}\,Q(k_1)\cap{\rm Dom}\,P(\overline{k}_2)\cap{\rm Dom}\,Q(\overline{k}_1)$ be a unit vector. We denote by $\mathbb{P}_\psi$ the orthogonal projection onto the one-dimensional subspace spanned by $\psi$. Evaluating the commutator $[Q(k_1)-P(k_2),\mathbb{P}_\psi]$ on a vector $\phi\in{\rm Dom}\,P(k_2)\cap{\rm Dom}\,Q(k_1)$, we get
\begin{eqnarray*}
[Q(k_1)-P(k_2),\mathbb{P}_\psi]\phi &=& \langle\psi,\phi\rangle \big(Q(k_1)-P(k_2)\big)(\psi) - \langle\psi,\big(Q(k_1)-P(k_2)\big)(\phi)\rangle \psi \\
&=& \langle\psi,\phi\rangle \big(Q(k_1)-P(k_2)\big)(\psi) - \langle\big(Q(\overline{k}_1)-P(\overline{k}_2)\big)\psi,\phi\rangle \psi
\end{eqnarray*}
We see that the range of $[Q(k_1)-P(k_2),\mathbb{P}_\psi]$ is two-dimensional, so it can be extended to a bounded operator on $\eufrak{H}$. Therefore $\mathbb{P}_\psi\in {\rm Dom}\, D_k$, and we get
\[
(D_k\mathbb{P}_\psi)\phi =\frac{i}{2}\Big(\langle\psi,\phi\rangle \big(Q(k_1)-P(k_2)\big)(\psi) - \langle\big(Q(\overline{k}_1)-P(\overline{k}_2)\big)\psi,\phi\rangle \psi\Big)
\]
for all $\phi\in\eufrak{H}$.
\end{example}

\begin{example}
Let $h\in \eufrak{h}\otimes\mathbb{R}^2$, $k\in \eufrak{h}_\mathbb{C}\otimes\mathbb{C}^2$. Then $\frac{i}{2}[Q(k_1)-P(k_2),U(h_1,h_2)]$ defines a bounded operator on $\eufrak{H}$, and we get
\[
D_k U(h_1,h_2) = i \big(\langle \overline{k}_1, h_1\rangle+\langle \overline{k}_2,h_2 \rangle\big) U(h_1,h_2).
\]
\end{example}

\begin{proposition}\label{prop closable1}
Let $k\in \eufrak{h}_\mathbb{C}\otimes\mathbb{C}^2$. The operator $D_k$ is a closable operator from $\mathcal{B}(\eufrak{H})$ to $\mathcal{B}(\eufrak{H})$ with respect to the weak topology.
\end{proposition}
\begin{proof}
Let $(B_n)_{n\in\mathbb{N}}\subseteq{\rm Dom}\,D_k\subseteq\mathcal{B}(\eufrak{H})$ be any sequence such that $B_n\to 0$ and $D_kB_n\to \beta$ for some $\beta\in\mathcal{B}(\eufrak{H})$ in the weak topology. To show that $D_k$ is closable, we have to show that this implies $\beta=0$. Let us evaluate $\beta$ between two exponential vectors $\mathcal{E}(h_1)$, $\mathcal{E}(h_2)$, $h_1,h_2\in\eufrak{h}_\mathbb{C}$, then we get
\begin{eqnarray*}
\langle \mathcal{E}(h_1),\beta \mathcal{E}(h_2)\rangle &=& \lim_{n\to \infty} \langle \mathcal{E}(h_1), D_kB_n\mathcal{E}(h_2)\rangle \\
&=& \lim_{n\to \infty}  \frac{i}{2} \big\langle \big(Q(\overline{k}_1)-P(\overline{k}_2)\big)\mathcal{E}(h_1), B_n\mathcal{E}(h_2)\big\rangle \\
&&-\lim_{n\to\infty}  \frac{i}{2}\big\langle\mathcal{E}(h_1), B_n\big(Q(k_1)-P(k_2)\big)\mathcal{E}(h_2)\big\rangle \\
&=& 0,
\end{eqnarray*}
and therefore $\beta=0$, as desired. 
\end{proof}

\begin{definition}\label{def derivation2}
We set
\[
\mathcal{S}= {\rm alg}\,\left\{O_h(\varphi)\Big| h\in \eufrak{h}\otimes\mathbb{R}; \varphi\in C^\infty(\mathbb{R}^2) \text{ s.t. }\frac{\partial^{\kappa_1+\kappa_2} \varphi}{\partial x^{\kappa_1}\partial y^{\kappa_2}}\in {\rm Dom}\, O_h \text{ for all } \kappa_1,\kappa_2\ge 0\right\},
\]
the elements of $\mathcal{S}$ will play the role of the smooth functionals. Note that $\mathcal{S}$ is weakly dense in $\mathcal{B}(\eufrak{H})$, i.e.\ $\mathcal{S}''=B(\eufrak{H})$, since $\mathcal{S}$ contains the Weyl operators $U(h_1,h_2)$ with $h_1,h_2\in\eufrak{h}$.

We define $D:\mathcal{S}\to\mathcal{B}(\eufrak{H})\otimes \eufrak{h}_\mathbb{C}\otimes \mathbb{C}^2$ (where the tensor product is the algebraic tensor product over $\mathbb{C}$) by setting $DO_h(\varphi)$ equal to
\[
DO_h(\varphi)= \left(\begin{array}{c} O_h\left(\displaystyle\frac{\partial \varphi}{\partial x}\right)\otimes h_1 \\[4mm] O_h\left(\displaystyle\frac{\partial \varphi}{\partial y}\right)\otimes h_2 \end{array}\right)
\]
and extending it as a derivation w.r.t.\ the $\mathcal{B}(\eufrak{H})$-bimodule structure of $\mathcal{B}(\eufrak{H})\otimes\eufrak{h}_\mathbb{C}\otimes \mathbb{C}^2$ defined by
\[
O\cdot \left(\begin{array}{c} O_1\otimes k_1 \\ O_2\otimes k_2 \end{array}\right) = \left(\begin{array}{c} OO_1\otimes k_1 \\ OO_2\otimes k_2 \end{array}\right), \quad
\left(\begin{array}{c} O_1\otimes k_1 \\ O_2\otimes k_2 \end{array}\right)\cdot O = \left(\begin{array}{c} O_1O\otimes k_1 \\ O_2O\otimes k_2 \end{array}\right) \
\]
for $O, O_1,O_2\in\mathcal{B}(\eufrak{H})$ and $k\in\eufrak{h}_\mathbb{C}\otimes\mathbb{C}^2$.
\end{definition}

\begin{example}
For $h\in \eufrak{h}\otimes\mathbb{R}^2$, we get
\begin{eqnarray*}
D U(h_1,h_2)&=& D O_h\Big(\exp i(x+y)\Big)
= i\left(\begin{array}{c}U(h_1,h_2)\otimes h_1 \\ U(h_1,h_2)\otimes h_2\end{array}\right) \\
&=& i U(h_1,h_2)\otimes h.
\end{eqnarray*}
\end{example}

\begin{definition}
We can define a $B(\eufrak{H})$-valued inner product on $\mathcal{B}(\eufrak{H})\otimes\eufrak{h}_\mathbb{C}\otimes\mathbb{C}^2$ by $\langle\cdot,\cdot\rangle:\mathcal{B}(\eufrak{H})\otimes\eufrak{h}_\mathbb{C}\otimes\mathbb{C}^2\times \mathcal{B}(\eufrak{H})\otimes\eufrak{h}_\mathbb{C}\otimes\mathbb{C}^2 \to \mathcal{B}(\eufrak{H})$ by
\[
\left\langle \left(\begin{array}{c} O_1\otimes h_1 \\ O_2\otimes h_2\end{array}\right), \left(\begin{array}{c} O'_1\otimes k_1 \\ O'_2\otimes k_2\end{array}\right)\right\rangle =  O^*_1O'_1 \langle h_1, k_1\rangle + O^*_2O'_2 \langle h_2, k_2\rangle
\]
\end{definition}
We have
\begin{eqnarray*}
\langle B,A\rangle &=& \langle A,B\rangle^* \\
O^*\langle A,B\rangle &=& \langle AO,B\rangle \\
\langle A,B\rangle O &=& \langle A,BO\rangle \\
\langle O^*A,B\rangle &=& \langle A,OB\rangle
\end{eqnarray*}
for all $A,B\in \mathcal{B}(\eufrak{H})\otimes\eufrak{h}_\mathbb{C}\otimes\mathbb{C}^2$ and all $O\in \mathcal{B}(\eufrak{H})$. This turns $\mathcal{B}(\eufrak{H})\otimes\eufrak{h}_\mathbb{C}\otimes\mathbb{C}^2$ into a pre-Hilbert module over $\mathcal{B}(\eufrak{H})$. It can be embedded in the Hilbert module $\eufrak{M}=\mathcal{B}(\eufrak{H},\eufrak{H}\otimes\eufrak{h}_\mathbb{C}\otimes\mathbb{C}^2)$ by mapping $O\otimes k\in\mathcal{B}(\eufrak{H})\otimes\eufrak{h}_\mathbb{C}\otimes\mathbb{C}^2 $ to the linear map $\eufrak{H}\ni v \mapsto Ov\otimes k\in \eufrak{H}\otimes\eufrak{h}_\mathbb{C}\otimes\mathbb{C}^2$. We will regard $\eufrak{h}_\mathbb{C}\otimes \mathbb{C}^2$ as a subspace of $\eufrak{M}$ via the embedding $\eufrak{h}_\mathbb{C}\ni k\mapsto {\rm id}_\eufrak{H}\otimes k\in \eufrak{M}$. Note that we have $O\cdot k=k\cdot O = O\otimes k$ and $\langle A,k\rangle=\langle \overline{k}, \overline{A}\rangle$ for all $k\in\eufrak{h}_\mathbb{C}\otimes\mathbb{C}^2$, $O\in\mathcal{B}(\eufrak{H})$, $A\in \eufrak{M}$, where the conjugation in $\eufrak{M}$ is defined by $\overline{O\otimes k}=O^*\otimes\overline{k}$.

\begin{proposition}\label{prop frechet}
Let $O\in \mathcal{S}$ and $k\in \eufrak{h}_\mathbb{C}\otimes\mathbb{C}^2$. Then $O\in {\rm Dom}\,D_k$ and
\[
D_k O = \langle \overline{k}, DO\rangle = \langle \overline{DO},k \rangle .
\]
\end{proposition}
\begin{proof}
For $h\in \eufrak{h}\otimes\mathbb{R}^2$ and $\varphi\in {\rm Dom}\,O_h$ s.t.\ also $\frac{\partial\varphi}{\partial x},\frac{\partial\varphi}{\partial y}\in {\rm Dom}\,O_h$, we get
\begin{eqnarray*}
\langle \overline{k}, DO_h(\varphi)\rangle &=& \left\langle \left(\begin{array}{c} \overline{k}_1 \\ \overline{k}_2\end{array}\right), \left(\begin{array}{c} O_h\left(\displaystyle\frac{\partial \varphi}{\partial x}\right)\otimes h_1 \\[4mm] O_h\left(\displaystyle\frac{\partial \varphi}{\partial y}\right)\otimes h_2\end{array}\right)\right\rangle \\
&=& O_h\left(\langle \overline{k}_1, h_1\rangle\frac{\partial\varphi}{\partial x}+ \langle \overline{k}_2, h_2\rangle\frac{\partial\varphi}{\partial y}\right)
\\
&=& \frac{i}{2}[Q(k_1)-P(k_2),O_h(\varphi)] = D_k O,
\end{eqnarray*}
where we used Proposition \ref{prop-int by parts}. The first equality of the proposition now follows, since both $O\mapsto D_kO=\frac{i}{2}[Q(k_1)-P(k_2),O]$ and $O\mapsto \langle \overline{k}, DO\rangle$ are derivations.

The second equality follows immediately.
\end{proof}

The next result is the analogue of Equation \eqref{int by parts1}.
\begin{theorem}\label{theo int by parts}
We have
\[
\mathbb{E}\big( \langle \overline{k}, DO\rangle \big) = \frac{1}{2}\mathbb{E}\big(\{P(k_1)+Q(k_2),O\}\big)
\]
for all $k\in \eufrak{h}_\mathbb{C}\otimes\mathbb{C}^2$ and all $O\in \mathcal{S}$, where $\{\cdot,\cdot\}$ denotes the anti-commutator $\{X,Y\}=XY+YX$.
\end{theorem}
\begin{proof}
This formula is a consequence of the fact that $Q(h)\Omega=h=iP(h)\Omega$ for all $h\in\eufrak{h}_\mathbb{C}$, we get
\begin{eqnarray*}
\mathbb{E}\big( \langle \overline{k}, DO\rangle \big) &=& \frac{i}{2} \Big( \langle \big(Q(\overline{k}_1)-P(\overline{k}_2)\big)\Omega,O\Omega\rangle - \langle\Omega, O\big(Q(k_1)-P(k_2)\big)\Omega\rangle\Big) \\
&=&  \frac{i}{2}\big(\langle \overline{k}_1+i\overline{k}_2 , O\Omega\rangle - \langle \Omega, O (k_1+ik_2)\rangle\big) \\
 &=& \frac{1}{2}\Big( \langle \big(P(\overline{k}_1)+Q(\overline{k}_2)\big)\Omega,O\Omega\rangle + \langle\Omega, O\big(P(k_1)+Q(k_2)\big)\Omega\rangle\Big) \\
&=&\frac{1}{2}\mathbb{E}\big(\{P(k_1)+Q(k_2),O\}\big).
\end{eqnarray*}
\end{proof}

There is also an analogue of \eqref{int by parts2}.
\begin{corollary}\label{cor int by parts}
Let $k\in \eufrak{h}_\mathbb{C}\otimes\mathbb{C}^2$, and $O_1,\ldots, O_n\in\mathcal{S}$, then
\[
\frac{1}{2}\mathbb{E}\left(\left\{P(k_1)+Q(k_2),\prod_{m=1}^n O_m\right\}\right) = \mathbb{E}\left( \sum_{m=1}^n \prod_{j=1}^{m-1} O_j\langle \overline{k}, DO_m\rangle \prod_{j=m+1}^{n} O_j\right),
\]
where the products are ordered such that the indices increase from the left to the right.
\end{corollary}
\begin{proof}
This is obvious, since $O\mapsto \langle \overline{k}, DO\rangle$ is a derivation.
\end{proof}

This formula for $n=3$ can be used to show that $D$ is a closable operator from $\mathcal{B}(\eufrak{H})$ to $\eufrak{M}$.
\begin{corollary}\label{cor closable2}
The derivation operator $D$ is a closable operator from $\mathcal{B}(\eufrak{H})$ to the $\mathcal{B}(\eufrak{H})$-Hilbert module $\eufrak{M}=\mathcal{B}(\eufrak{H},\eufrak{H}\otimes \eufrak{h}_\mathbb{C}\otimes\mathbb{C}^2)$ w.r.t.\ the weak topologies.
\end{corollary}
\begin{proof}
We have to show that for any sequence $(A_n)_{n\in\mathbb{N}}$ in $\mathcal{S}$ with $A_n\to 0$ and $DA_n\to \alpha \in \eufrak{M}$, we get $\alpha=0$. Let $f,g\in\eufrak{h}_\mathbb{C}$. Set $f_1=\frac{f+\overline{f}}{2}$, $f_2=\frac{f-\overline{f}}{2i}$, $g_1=\frac{g+\overline{g}}{2}$, and $g_2=\frac{g-\overline{g}}{2i}$, then we have $U(f_1,f_2)\Omega=e^{-||f||/2}\mathcal{E}(f)$ and $U(g_1,g_2)\Omega=e^{-||g||/2}\mathcal{E}(g)$. Thus we get
\begin{eqnarray*}
&&\exp\big((||f||^2+||g||^2)/2\big)\langle \mathcal{E}(f)\otimes \overline{h},\alpha \mathcal{E}(g)\rangle \\
&=& \exp\big((||f||^2+||g||^2)/2\big)\langle \mathcal{E}(f),\langle \overline{h},\alpha\rangle \mathcal{E}(g)\rangle \\
&=& \lim_{n\to\infty}\mathbb{E} \big(U(-f_1,-f_2) \langle \overline{h},D A_n\rangle U(g_1,g_2)\big) \\
&=& \lim_{n\to\infty}\mathbb{E}\Big( \frac{1}{2}\big\{P(h_1)+Q(h_2),U(-f_1,-f_2)A_n U(g_1,g_2)\big\} \\
&& - \big\langle \overline{h}, DU(-f_1,-f_2)\big\rangle A_n U(g_1,g_2)
- U(-f_1,-f_2)A_n\big\langle \overline{h}, DU(g_1,g_2)\big\rangle\Big) \\
&=& \lim_{n\to\infty}\big( \langle \psi_1, A_n \psi_2\rangle +\langle \psi_3, A_n \psi_4\rangle -\langle \psi_5, A_n \psi_6\rangle -\langle \psi_7, A_n \psi_8\rangle\big)  \\
&=& 0
\end{eqnarray*}
for all $h\in\eufrak{h}_\mathbb{C}\otimes\mathbb{C}^2$, where
\[
\begin{array}{rclcrcl}
\psi_1&=&\frac{1}{2}U(f_1,f_2)\big(P(\overline{h}_1)+Q(\overline{h}_2)\big)\Omega, && \psi_2&=& U(g_1,g_2)\Omega, \\
\psi_3&=&U(f_1,f_2)\Omega, && \psi_4&=& \frac{1}{2}U(g_1,g_2)\big(P(h_1)+Q(h_2)\big)\Omega, \\
\psi_5&=&\big(D_hU(-f_1,-f_2)\big)^*\Omega, && \psi_6&=&  U(g_1,g_2)\Omega, \\
\psi_7&=&U(f_1,f_2)\Omega, && \psi_8&=&D_hU(g_1,g_2)\Omega.
\end{array}
\]
But this implies $\alpha=0$, since $\{\mathcal{E}(f)\otimes \overline{h}|f\in \eufrak{h}_\mathbb{C},h\in\eufrak{h}_\mathbb{C}\otimes\mathbb{C}^2\}$ is dense in $\eufrak{H}\otimes\eufrak{h}_\mathbb{C}\otimes\mathbb{C}^2$.
\end{proof}

\begin{remark}
This implies that $D$ is also closable in stronger topologies, such as, e.g., the norm topology and the strong topology.
\end{remark}

We will denote the closure of $D$ again by the same symbol.

\begin{proposition}
Let $O\in{\rm Dom}\,D$. Then $O^*\in{\rm Dom}\,D$ and
\[
DO^* = \overline{DO}.
\]
In particular, since $D$ is a derivation, this implies that ${\rm Dom}\,D$ is a $*$-subalgebra of $\mathcal{B}(\eufrak{H})$.
\end{proposition}
\begin{proof}
It is not difficult to check this directly on the Weyl operators $U(h_1,h_2)$, $h\in\eufrak{h}\otimes\mathbb{R}^2$. We get $U(h_1,h_2)^*=U(-h_1,-h_2)$ and
\begin{eqnarray*}
D\big(U(h_1,h_2)^*\big) &=& DU(-h_1,-h_2)= -i U(-h_1,-h_2)\otimes h \\
&=& U(h_1,h_2)^* \otimes \overline{(ih)} = \overline{DU(h_1,h_2)}.
\end{eqnarray*}
By linearity and continuity it therefore extends to all of ${\rm Dom}\,D$.
\end{proof}

We will now show how $D$ can be iterated. Let $H$ be a complex Hilbert space, then we can define a derivation operator $D:\mathcal{S}\otimes H\to \mathcal{B}(\eufrak{H})\otimes\eufrak{h}_\mathbb{C}\otimes\mathbb{C}^2\otimes H$ by setting $D(O\otimes h)=DO\otimes h$ for $O\in\mathcal{S}$ and $h\in H$. Closing it, we get an unbounded derivation from the Hilbert module $\mathcal{B}(\eufrak{H}, \eufrak{H}\otimes H)$ to $\eufrak{M}(H)=\mathcal{B}(\eufrak{H}\otimes H,\eufrak{H}\otimes\eufrak{h}_\mathbb{C}\otimes \mathbb{C}^2\otimes H)$. This allows us to iterate $D$. It is easy to see that $D$ maps $\mathcal{S}\otimes H$ to $\mathcal{S}\otimes\eufrak{h}_\mathbb{C}\otimes \mathbb{C}^2\otimes H$ and so we have $D^n(\mathcal{S}\otimes H)\subseteq \mathcal{S}\otimes \left(\eufrak{h}_\mathbb{C}\otimes \mathbb{C}^2\right)^{\otimes n}\otimes H$. In particular, $\mathcal{S}\subseteq{\rm Dom}\,D^n$ for all $n\in\mathbb{N}$, and we can define Sobolev-type norms $||\cdot||_n$ and semi-norms $||\cdot||_{\psi,n}$, on $\mathcal{S}$ by
\begin{eqnarray*}
||O||^2_n &=& ||O^*O||+\sum_{j=1}^n ||\langle D^n O,D^nO\rangle||, \\
||O||^2_{\psi,n} &=& ||O\psi||^2+\sum_{j=1}^n ||\langle \psi, \langle D^n O,D^nO\rangle\psi\rangle||, \quad\text{ for } \psi\in\eufrak{H}
\end{eqnarray*}
In this way we can define Sobolev-type topologies on ${\rm Dom}\,D^n$.

We will now extend the definition of the ``Fr\'echet derivation'' $D_k$ to the case where $k$ is replaced by an element of $\eufrak{M}$. It becomes now important to distinguish between a right and a left ``derivation operator''. Furthermore, it is no longer a derivation.
\begin{definition}
Let $u\in\eufrak{M}$ and $O\in {\rm Dom}\,D$. Then we define the right gradient $\overrightarrow{D}_uO$ and the left gradient $O\overleftarrow{D}_u$ of $O$ with respect to $u$ by
\begin{eqnarray*}
\overrightarrow{D}_u O&=& \langle \overline{u},DO\rangle, \\
O\overleftarrow{D}_u &=& \langle \overline{DO},u\rangle. \\
\end{eqnarray*}
\end{definition}

We list several properties of the gradient.
\begin{proposition}
\begin{enumerate}
\item
Let $X\in\mathcal{B}(\eufrak{H})$, $O,O_1,O_2\in {\rm Dom}\,D$, and $u\in\eufrak{M}$. Then
\begin{eqnarray*}
\overrightarrow{D}_{Xu}O &=& X \overrightarrow{D}_uO, \\
\overrightarrow{D}_u(O_1O_2) &=& \left(\overrightarrow{D}_uO_1\right) O_2+\overrightarrow{D}_{uO_1}O_2 , \\
O\overleftarrow{D}_{uX} &=& (O\overleftarrow{D}_u)X, \\
(O_1O_2)\overleftarrow{D}_u &=& O_1\overleftarrow{D}_{O_2u} + O_1\left(O_2\overleftarrow{D}_u\right) , \\
\end{eqnarray*}
\item
For $k\in \eufrak{h}_\mathbb{C}\otimes \mathbb{C}^2$ and $O\in{\rm Dom}\,D$, we have
\[
D_k O = \overrightarrow{D}_{{\rm id}_\eufrak{H}\otimes k}O = O\overleftarrow{D}_{{\rm id}_\eufrak{H}\otimes k}
\]
\end{enumerate}
\end{proposition}
\begin{proof}
These properties can be deduced easily from the definition of the gradient and the properties of the derivation operator $D$ and the inner product $\langle,\rangle$.
\end{proof}

\begin{remark}
We can also define a two-sided gradient $\overleftrightarrow{D}_u:{\rm Dom}\,D\times{\rm Dom}\,D\to \mathcal{B}(\eufrak{H})$ by $\overleftrightarrow{D}_u:(O_1,O_2)\mapsto O_1\overleftrightarrow{D}_uO_2= O_1\left(\overrightarrow{D}_uO_2\right) + \left(O_1\overleftarrow{D}_u\right)O_2$. For $k\in \eufrak{h}_\mathbb{C}\otimes \mathbb{C}^2$ we have $O_1\overleftrightarrow{D}_{{\rm id}_\eufrak{H}\otimes k}  O_2= D_k(O_1O_2)$.
\end{remark}

\section{The divergence operator}

The algebra $\mathcal{B}(\eufrak{H})$ of bounded operators on the symmetric Fock space $\eufrak{H}$ and the Hilbert module $\eufrak{M}$ are not Hilbert spaces with respect to the expectation in the vacuum vector $\Omega$. Therefore we can not define the divergence operator or Skorohod integral $\delta$ as the adjoint of the derivation $D$. It might be tempting to try to define $\delta X$ as an operator such that the condition
\begin{equation}\label{duality}
\mathbb{E}\big((\delta X) B\big)\overset{?}{=} \mathbb{E}\left(\overrightarrow{D}_X B\right)
\end{equation}
is satisfied for all $B\in{\rm Dom}\,\overrightarrow{D}_X$, even though it is not sufficient to characterize $\delta X$. But the following proposition shows that this is not possible.

\begin{proposition}\label{prop no go}
Let $k\in \eufrak{h}_\mathbb{C}\otimes\mathbb{C}^2$ with $k_1+ik_2\not=0$. There exists no (possibly unbounded) operator $M$ whose domain contains the vacuum vector such that
\[
\mathbb{E}\big( M B\big) = \mathbb{E}\left(D_k B\right)
\]
holds for all $B\in{\rm Dom}\,D_k$.
\end{proposition}
\begin{proof}
We assume that such an operator $M$ exists and show that this leads to a contradiction.

Let $B\in\mathcal{B}(\eufrak{H})$ be the operator defined by $\eufrak{H}\ni\psi\mapsto \langle k_1+ik_2,\psi\rangle\Omega$, it is easy to see that $B\in{\rm Dom}\,D_k$ and that $D_kB$ is given by
\[
(D_kB)\psi =  \frac{i}{2}\langle k_1+ik_2,\psi\rangle(k_1+ik_2)-\frac{i}{2}\big\langle\big(Q(\overline{k}_1)-P(\overline{k}_2)\big)(k_1+ik_2),\psi\big\rangle \Omega,\text{ for }\psi\in\eufrak{H}.
\]
Therefore, if $M$ existed, we would have
\begin{eqnarray*}
0 &=& \langle \Omega, MB\Omega\rangle = \mathbb{E}(MB) = \mathbb{E}(D_kB) \\
&=& \langle \Omega,(D_kB)\Omega\rangle =  -\frac{i}{2}\langle k_1+ik_2,k_1+ik_2\rangle,
\end{eqnarray*}
which is clearly impossible.
\end{proof}

We introduce the analogue of smooth elementary $\eufrak{h}$-valued random variables,
\[
\mathcal{S}_\eufrak{h}=\left\{ u=\sum_{j=1}^n F_j\otimes h^{(j)}\Big| n\in\mathbb{N}, F_1,\ldots,F_n\in\mathcal{S},h^{(1)},\ldots,h^{(n)}\in\eufrak{h}_\mathbb{C}\otimes\mathbb{C}^2\right\}.
\]
If we define $A\overleftrightarrow{\delta u}B$ for $u=\sum_{j=1}^n F_j\otimes h^{(j)}\in\mathcal{S}_\eufrak{h}$ and $A,B\in B(\eufrak{H})$ by
\[
A\overleftrightarrow{\delta u}B= \frac{1}{2}\sum_{j=1}^n \left\{ P\left(h^{(j)}_1\right)+Q\left(h^{(j)}_2\right) ,AF_jB\right\} - A\sum_{j=1}^n \left(D_{h^{(j)}}F_j\right)B.
\]
then it follows from Corollary \ref{cor int by parts} that this satisfies
\begin{equation}\label{duality2}
\mathbb{E}\left(A\overleftrightarrow{\delta u}B\right) = \mathbb{E}\left(A\overleftrightarrow{D}_u B\right).
\end{equation}
But this can only be written as a product $AXB$ with some operator $X$, if $A$ and $B$ commute with $P\left(h^{(1)}_1\right)+Q\left(h^{(1)}_2\right),\ldots, P\left(h^{(n)}_1\right)+Q\left(h^{(n)}_2\right)$. We see that a condition of the form (\ref{duality}) or (\ref{duality2}) is too strong, if we require it to be satisfied for all $A,B\in{\rm Dom}\,D$. We have to impose some commutativity on $A$ and $B$ to weaken the condition, in order to be able to satisfy it. We will now give a first definition of a divergence operator that satisfies a weaker version of (\ref{duality2}), see Proposition \ref{prop duality} below. In Remark \ref{ref def delta} we will extend this definition to a bigger domain.

\begin{definition}\label{def skorohod1}
We set 
\begin{eqnarray*}
\mathcal{S}_{\eufrak{h},\delta} &=& \Big\{ u=\sum_{j=1}^n F_j\otimes h^{(j)}\in\mathcal{S}_\eufrak{h}\Big| \frac{1}{2}\sum_{j=1}^n \left\{ P\left(h^{(j)}_1\right)+Q\left(h^{(j)}_2\right) ,F_j\right\} - \sum_{j=1}^n D_{h^{(j)}}F_j\\
&& \qquad\qquad\qquad\qquad\qquad \text{ defines a bounded operator on }\eufrak{H}\Big\}
\end{eqnarray*}
and define the divergence operator $\delta:\mathcal{S}_{\eufrak{h},\delta}\to B(\eufrak{H})$ by 
\[
\delta(u) = \frac{1}{2}\sum_{j=1}^n \left\{ P\left(h^{(j)}_1\right)+Q\left(h^{(j)}_2\right) ,F_j\right\} - \sum_{j=1}^n D_{h^{(j)}}F_j.
\]
for $u=\sum_{j=1}^n F_j\otimes h^{(j)}\in\mathcal{S}_{\eufrak{h},\delta}$.
\end{definition}

It is easy to check that
\[
\delta(\overline{u}) = \big(\delta(u)\big)^*
\]
holds for all $u\in\mathcal{S}_{\eufrak{h},\delta}$.

\begin{proposition}\label{prop duality}
Let $u=\sum_{j=1}^n F_j\otimes h^{(j)}\in\mathcal{S}_{\eufrak{h},\delta}$ and
\[
A,B\in {\rm Dom}\,D\cap\left\{ P\left(h^{(1)}_1\right)+Q\left(h^{(1)}_2\right),\ldots, P\left(h^{(n)}_1\right)+Q\left(h^{(n)}_2\right)\right\}'
\]
i.e., $A$ and $B$ are in the commutant of $\left\{ P\left(h^{(1)}_1\right)+Q\left(h^{(1)}_2\right),\ldots, P\left(h^{(n)}_1\right)+Q\left(h^{(n)}_2\right)\right\}$, then we have
\[
\mathbb{E}\big(A\delta(u)B\big) = \mathbb{E}\left(A\overleftrightarrow{D}_u B\right).
\]
\end{proposition}
\begin{remark}
Note that $\delta:\mathcal{S}_\eufrak{h,\delta}\to\mathcal{B}(\eufrak{H})$ is the only linear map with this property, since for one single element $h\in\eufrak{h}_\mathbb{C}\otimes\mathbb{C}^2$, the sets
\begin{eqnarray*}
&\left\{ A^*\Omega \Big| A\in {\rm Dom}\,D\cap\left\{ P\left(h_1\right)+Q\left(h_2\right)\right\}'\right\}& \\
&\left\{ B\Omega \Big| B\in {\rm Dom}\,D\cap\left\{ P\left(h_1\right)+Q\left(h_2\right)\right\}'\right\}&
\end{eqnarray*}
are still total in $\eufrak{H}$.
\end{remark}
\begin{proof}
{}From Corollary \ref{cor int by parts} we get
\begin{eqnarray*}
\mathbb{E}\left(A\overleftrightarrow{D}_u B\right) &=&\mathbb{E}\big(A\langle \overline{u}, DB\rangle + \langle \overline{DA}, u\rangle B \big) \\
&=&\mathbb{E}\left( \sum_{j=1}^n AF_j \left( D_{h^{(j)}}B\right) +  \sum_{j=1}^n \left( D_{h^{(j)}}A\right)F_jB\right) \\
&=& \mathbb{E} \left(\frac{1}{2}\sum_{j=1}^n \left\{ P\left(h^{(j)}_1\right)+Q\left(h^{(j)}_2\right) ,AF_jB\right\} - \sum_{j=1}^n A \left(D_{h^{(j)}}F_j\right)B\right).
\end{eqnarray*}
But since $A$ and $B$ commute with $P\left(h^{(1)}_1\right)+Q\left(h^{(1)}_2\right),\ldots, P\left(h^{(n)}_1\right)+Q\left(h^{(n)}_2\right)$, we can pull them out of the anti-commutator, and we get
\begin{eqnarray*}
\mathbb{E}\left(A\overleftrightarrow{D}_u B\right) &=&\mathbb{E} \left(\frac{1}{2}\sum_{j=1}^n A \left\{ P\left(h^{(j)}_1\right)+Q\left(h^{(j)}_2\right) ,F_j\right\}B - \sum_{j=1}^n A\left(D_{h^{(j)}}F_j\right)B\right) \\
&=& \mathbb{E}(A\delta(u)B).
\end{eqnarray*}
\end{proof}

We will now give an explicit formula for the matrix elements between two exponential vectors of the divergence of a smooth elementary element $u\in \mathcal{S}_{\eufrak{h},\delta}$, this is the analogue of the first fundamental lemma in the Hudson-Parthasarathy calculus, see, e.g., \cite[Proposition 25.1]{parthasarathy92}.

\begin{theorem}\label{theo exp values}
Let $u\in\mathcal{S}_{\eufrak{h},\delta}$. Then we have the following formula
\[
\langle\mathcal{E}(k_1), \delta(u)\mathcal{E}(k_2)\rangle = \left\langle\mathcal{E}(k_1)\otimes\left(\begin{array}{c} ik_1-i\overline{k}_2 \\ k_1+\overline{k}_2\end{array}\right)   , u \mathcal{E}(k_2)\right\rangle
\]
for the evaluation of the divergence $\delta(u)$ of $u$ between two exponential vectors $\mathcal{E}(k_1)$, $\mathcal{E}(k_2)$, for $k_1,k_2\in\eufrak{h}_\mathbb{C}$.
\end{theorem}
\begin{remark}\label{ref def delta}
This suggests to extend the definition of $\delta$ in the following way: set
\begin{eqnarray}
{\rm Dom}\,\delta &=& \Big\{ u\in \eufrak{M}\Big| \exists M\in\eufrak{B}(\eufrak{H}), \forall k_1,k_2\in\eufrak{h}_\mathbb{C}: \nonumber\\
&& \qquad\qquad \langle\mathcal{E}(k_1), M \mathcal{E}(k_2)\rangle = \left\langle\mathcal{E}(k_1)\otimes\left(\begin{array}{c} ik_1-i\overline{k}_2 \\ k_1+\overline{k}_2\end{array}\right)   , u \mathcal{E}(k_2)\right\rangle\Big\} \label{eq-def-delta}
\end{eqnarray}
and define $\delta(u)$ for $u\in{\rm Dom}\,\delta$ to be the unique operator $M$ that satisfies the condition in Equation (\ref{eq-def-delta}).
\end{remark}
\begin{proof}
Let $u=\sum_{j=1}^n F_j\otimes h^{(j)}$. Recalling the definition of $D_h$ we get the following alternative expression for $\delta(u)$,
\begin{eqnarray}
\delta(u) &=&\frac{1}{2} \sum_{j=1}^n \Big(P\left(h^{(j)}_1\right)+Q\left(h^{(j)}_2\right)-iQ\left(h^{(j)}_1\right)+iP\left(h^{(j)}_2\right)\Big) F_j \nonumber \\
&& +\frac{1}{2}\sum_{j=1}^n F_j \Big(P\left(h^{(j)}_1\right) +Q\left(h^{(j)}_2\right)+iQ\left(h^{(j)}_1\right)-iP\left(h^{(j)}_2\right)\Big) \nonumber \\
&=&\sum_{j=1}^n \Big( a^+\left(h^{(j)}_2-ih^{(j)}_1\right)F_j + F_j a\left(\overline{h^{(j)}_2-ih^{(j)}_1}\right)\Big). \label{wick}
\end{eqnarray}
Evaluating this between two exponential vectors, we obtain
\begin{eqnarray*}
\langle\mathcal{E}(k_1), \delta(u)\mathcal{E}(k_2)\rangle &= & \sum_{j=1}^n\big\langle a(h_2^{(j)}-ih_1^{(j)})\mathcal{E}(k_1), F_j\mathcal{E}(k_2)\big\rangle \\
&& + \sum_{j=1}^n \big\langle \mathcal{E}(k_1), F_ja(\overline{h_2^{(j)}+ih_1^{(j)}})\mathcal{E}(k_2)\big\rangle \\
&=&  \sum_{j=1}^n \big(\overline{\langle h_2^{(j)}-ih_1^{(j)},k_1\rangle} + \langle \overline{h_2^{(j)}+ih_1^{(j)}},k_2\rangle\big)\langle \mathcal{E}(k_1), F_j\mathcal{E}(k_2)\rangle \\
&=& \sum_{j=1}^n \big(\langle k_1,h_2^{(j)}-ih_1^{(j)}\rangle + \langle \overline{k}_2, h_2^{(j)}+ih_1^{(j)}\rangle\big)\langle \mathcal{E}(k_1), F_j\mathcal{E}(k_2)\rangle \\
&=& \left\langle\mathcal{E}(k_1)\otimes\left(\begin{array}{c} ik_1-i\overline{k}_2 \\ k_1+\overline{k}_2\end{array}\right)   , u \mathcal{E}(k_2)\right\rangle
.
\end{eqnarray*}
\end{proof}

\begin{corollary}
The divergence operator $\delta$ is closable in the weak topology.
\end{corollary}
\begin{proof}
Let $(u_n)_{n\in\mathbb{N}}$ be a sequence such that $u_n\to 0$ and $\delta(u_n)\to \beta\in\mathcal{B}(\eufrak{H})$ in the weak topology. Then we get
\begin{eqnarray*}
\langle\mathcal{E}(k_1), \beta \mathcal{E}(k_2)\rangle &=& \lim_{n\to\infty}\langle\mathcal{E}(k_1), \delta(u_n) \mathcal{E}(k_2)\rangle \\
&=& \lim_{n\to\infty}
 \left\langle\mathcal{E}(k_1)\otimes\left(\begin{array}{c} ik_1-i\overline{k}_2 \\ k_1+\overline{k}_2\end{array}\right)   , u_n \mathcal{E}(k_2)\right\rangle \\
&=& 0
\end{eqnarray*}
for all $k_1,k_2\in\eufrak{h}_\mathbb{C}$, and thus $\beta=0$.
\end{proof}

We have the following analogues of Equations (\ref{comm rel1}) and (\ref{DFu1}).

\begin{proposition}
Let $u,v\in \mathcal{S}_{\eufrak{h},\delta}$, $F\in\mathcal{S}$, $h\in\eufrak{h}_\mathbb{C}\otimes\mathbb{C}^2$, then we have
\begin{eqnarray}
D_h\circ \delta(u) &=& \langle \overline{h},u\rangle + \delta \circ D_h(u) \label{comm rel}\\
\delta(Fu) &=& F\delta(u) - F\overleftarrow{D}_u  +\frac{1}{2}  \sum_{j=1}^n \left[P\left(h^{(j)}_1\right)+Q\left(h^{(j)}_2\right),F\right] F_j  \label{DFu}\\
\delta(uF) &=&  \delta(u)F - \overrightarrow{D}_uF  + \frac{1}{2} \sum_{j=1}^n \left[F,P\left(h^{(j)}_1\right)+Q\left(h^{(j)}_2\right)\right] F_j  \label{DuF}
\end{eqnarray}
\end{proposition}
\begin{proof}
\begin{enumerate}
\item
Let $u=\sum_{j=1}^n F_j\otimes h^{(j)}$. We set
\begin{eqnarray*}
X_j &=& \frac{1}{2}\left(P\left(h^{(j)}_1\right)+Q\left(h^{(j)}_2\right)\right),\\
Y_j &=& \frac{i}{2}\left(Q\left(h^{(j)}_1\right)-P\left(h^{(j)}_2\right)\right),\\
 Y &=& \frac{i}{2}\big(Q\left(h_1\right)-P\left(h_2\right)\big),
\end{eqnarray*}
then we have
$\delta(u)=\sum_{j=1}^n \big((X_j-Y_j)F_j + F_j(X_j+Y_j)\big)$,
and therefore
\[
D_h\big(\delta(u)\big) = \sum_{j=1}^n \big(Y(X_j-Y_j)F_j + Y F_j(X_j+Y_j)- (X_j-Y_j)F_jY - F_j(X_j+Y_j)Y \big).
\]
On the other hand we have $D_h(u)=\sum_{j=1}^n (YF_j-F_jY)\otimes h^{(j)}$, and
\[
\delta\big(D_h(u)\big) = \sum_{j=1}^n \big((X_j-Y_j)YF_j- (X_j-Y_j)F_jY + YF_j(X_j+Y_j)- F_jY(X_j+Y_j)\big).
\]
Taking the difference of these two expressions, we get
\begin{eqnarray*}
D_h\big(\delta(u)\big)-\delta\big(D_h(u)\big) &=&\sum_{j=1}^n \big(\big[Y,X_j-Y_j\big]F_j+ F_j[Y,X_j+Y_j]\big) \\
&=& \sum_{j=1}^n \big(\langle \overline{h}_1,h^{(j)}_1\rangle + \langle\overline{h}_2,h^{(j)}_2\rangle\big) F_j = \langle \overline{h},u\rangle.
\end{eqnarray*}
\item
A straightforward computation gives
\begin{eqnarray*}
\delta(Fu)&=& \sum_{j=1}^n \big((X_j-Y_j)FF_j + FF_j(X_j+Y_j)\big) \\
&=&F \sum_{j=1}^n \big((X_j-Y_j)F_j + F_j(X_j+Y_j)\big) - \sum_{j=1}^n [F,X_j-Y_j] F_j \\
&=& F \delta(u) - \sum_{j=1}^n [Y_j,F] F_j +  \sum_{j=1}^n [X_jF,X] F_j \\
&=&   F \delta(u) - \sum_{j=1}^n \langle F^*,h^{(j)}\rangle F_j +  \sum_{j=1}^n [X_j,F] F_j \\
&=& F\delta(u) - F\overleftarrow{D}_u  +  \sum_{j=1}^n [X_j,F] F_j
\end{eqnarray*}
where we used that  $[X_j,F]=i\left\langle\left(\begin{array}{c}-\overline{h}_2 \\ \overline{h}_1\end{array}\right) ,DF\right\rangle$ defines a bounded operator, since $F\in\mathcal{S}\subseteq{\rm Dom}\,D$. Equation (\ref{DuF}) can be shown similarly.
\end{enumerate}
\end{proof}

If we impose additional commutativity conditions, which are always satisfied in the commutative case, then we get simpler formulas that are more similar to the classical ones.

\begin{corollary}
If $u=\sum_{j=1}^n F_j\otimes h^{(j)}\in\mathcal{S}_{\eufrak{h},\delta}$ and
\[
F\in {\rm Dom}\,D\cap\left\{ P\left(h^{(1)}_1\right)+Q\left(h^{(1)}_2\right),\ldots, P\left(h^{(n)}_1\right)+Q\left(h^{(n)}_2\right)\right\}'
\]
then we have
\begin{eqnarray*}
\delta(Fu) &=& F\delta(u) - F\overleftarrow{D}_u, \\
\delta(uF) &=&\delta(u)F - \overrightarrow{D}_uF.
\end{eqnarray*}
\end{corollary}

\section{Examples and applications}

\subsection{Relation to the commutative case}

In this section we will show that the non-commutative calculus studied here contains the commutative calculus, at least if we restrict ourselves to bounded functionals.

It is well-known that the symmetric Fock space $\Gamma(\eufrak{h}_\mathbb{C})$ is isomorphic to the complexification $L^2(\Omega;\mathbb{C})$ of the Wiener space $L^2(\Omega)$ over $\eufrak{h}$, cf.\ \cite{biane93,janson97,meyer95}. Such an isomorphism $I:L^2(\Omega;\mathbb{C})\overset{\cong}{\mapsto}\Gamma(\eufrak{h}_\mathbb{C})$ can be defined by extending the map
\[
I: e^{iW(h)}\mapsto I\left(e^{iW(h)}\right) = e^{iQ(h)}\Omega= e^{-||h||^2/2} \mathcal{E}(ih), \qquad \text{ for } h\in\eufrak{h}.
\]
Using this isomorphism, a bounded functional $F\in L^\infty(\Omega;\mathbb{C})$ becomes a bounded operator $M(F)$ on $\Gamma(\eufrak{h}_\mathbb{C})$, acting simply by multiplication,
\[
M(F)\psi= I\big(FI^{-1}(\psi)\big), \qquad \text{ for } \psi\in\Gamma(\eufrak{h}_\mathbb{C}).
\]
In particular, we get $M\left(e^{iW(h)}\right)=U(0,h)$ for $h\in\eufrak{h}$.

We can show that the derivation of a bounded differentiable functional coincides with its derivation as a bounded operator.
\begin{proposition}
Let $k\in\eufrak{h}$ and $F\in L^\infty(\Omega;\mathbb{C})\cap {\rm Dom}\,\tilde{D}_k$ s.t.\ $\tilde{D}_k F\in L^\infty(\Omega;\mathbb{C})$. Then we have $M(F)\in {\rm Dom}\,D_{k_0}$, where $k_0=\left(\begin{array}{c} 0 \\ k\end{array}\right)$, and
\[
M(\tilde{D}_k F) = D_{k_0} \big(M(F)\big).
\]
\end{proposition}
\begin{proof}
It is sufficient to check this for functionals of the form $F=e^{iW(h)}$, $h\in\eufrak{h}$. We get
\begin{eqnarray*}
M(\tilde{D}_k e^{iW(h)} ) &=& M\left(i\langle k,h\rangle e^{iW(h)}\right) \\
&=& i\langle k,h\rangle U(0,h) = i\left\langle\left(\begin{array}{c} 0 \\ k\end{array}\right),\left(\begin{array}{c} 0 \\ h\end{array}\right) \right\rangle U(0,h) \\
&=& D_{k_0} U(0,h) = D_ {k_0}\big(M(e^{iW(h)})\big).
\end{eqnarray*}
\end{proof}
This implies that we also have an analogous result for the divergence.

\subsection{Sufficient conditions for the existence of smooth densities}

In this section we will use the operator $D$ to give sufficient conditions for the existence and smoothness of densities for operators on $\eufrak{H}$. The first result is a generalisation of Proposition \ref{prop-regular1} to arbitrary states.

\begin{theorem}\label{theo regular2}
Let $\kappa\in\mathbb{N}$, $h\in\eufrak{h}\otimes\mathbb{R}^2$ with $\langle h_1,h_2\rangle\not=0$, and suppose that $\Phi$ is of the form 
\[
\Phi(X)={\rm tr}(\rho X) \text{ for all }X\in B(\eufrak{H}),
\]
for some density matrix $\rho$.
If there exist $k,\ell\in\eufrak{h}_\mathbb{C}\otimes\mathbb{C}^2$ such that
\[
\det\left(\begin{array}{cc}\langle h_1,k_1 \rangle & \langle h_2, k_2\rangle\\
\langle h_1,\ell_1 \rangle& \langle h_2, \ell_2\rangle\end{array}\right)\not=0,
\]
and $\rho\in\bigcap_{\kappa_1+\kappa_2\le\kappa} {\rm Dom}\,D_k^{\kappa_1}D_\ell^{\kappa_2}$, and
\[
{\rm tr}(|D_k^{\kappa_1}D_\ell^{\kappa_2}\rho|) < \infty \qquad\text{ for all } \kappa_1+\kappa_2\le \kappa,
\]
then $w_{h,\Phi}\in \bigcap_{2\le p\le\infty}H^{p,\kappa}(\mathbb{R}^2)$, i.e.\ the Wigner density $w_{h,\Phi}$ lies in the Sobolev spaces of order $\kappa$. 
\end{theorem}
\begin{proof}
Let
\[
A=\left(\begin{array}{cc} 
\langle h_1,k_1 \rangle & \langle h_2, k_2\rangle \\{}
\langle h_1,\ell_1 \rangle & \langle h_2, \ell_2\rangle
\end{array}\right)
\]
and set
\[
\left(\begin{array}{c} X_1 \\  X_2 \end{array}\right) = \frac{i}{2}\, A^{-1} \left(\begin{array}{c} Q(k_1)-P(k_2) \\ Q(\ell_1)-P(\ell_2) \end{array}\right), 
\]
then we have
\begin{eqnarray*}
[X_1,O_h(\varphi)]&=&\frac{\langle h_2, \ell_2\rangle D_k O_h(\varphi) - \langle h_2, k_2\rangle D_\ell O_h(\varphi)}{\det A}= O_h\left(\frac{\partial \varphi}{\partial x}\right), \\{}
[X_2,O_h(\varphi)]&=& \frac{-\langle h_1, \ell_1\rangle D_k O_h(\varphi) + \langle h_1, k_2\rangle D_\ell O_h(\varphi)}{\det A}=  O_h\left(\frac{\partial \varphi}{\partial y}\right),
\end{eqnarray*}
for all Schwartz functions $\varphi$. Therefore
\begin{eqnarray*}
\left|\int\frac{\partial^{\kappa_1+\kappa_2}\varphi}{\partial x^{\kappa_1}\partial y^{\kappa_2}}{\rm d}W_{h,\Phi}\right| &=& \left|{\rm tr}\left( \rho\, O_h\left( \frac{\partial^{\kappa_1+\kappa_2}\varphi}{\partial x^{\kappa_1}\partial y^{\kappa_2}}\right)\right)\right| \\
&=&\left|{\rm tr}\Big(\rho \underbrace{\big[ X_1,\ldots\big[X_1},\underbrace{\big[X_2,\ldots\big[X_2},O_h(\varphi)\big]\big]\big]\big]\Big)\right| \\
&&\hspace{17mm} \kappa_1 \text{ times }\hspace{7mm}\kappa_2\text{ times } \\[1mm]
&=&\left|{\rm tr}\big([ X_2,\ldots[X_2,[X_1,\ldots[X_1,\rho]]]]O_h(\varphi)\big)\right| \\[1mm]
&\le &C_{\rho,\kappa_1,\kappa_2}|| O_h(\varphi)|| \le C_{\rho,\kappa_1,\kappa_2} C_{h,p} ||\varphi||_p,
\end{eqnarray*}
for all $p\in[1,2]$, since $\rho\in\bigcap_{\kappa_1+\kappa_2\le\kappa} {\rm Dom}\,D_k^{\kappa_1}D_\ell^{\kappa_2}$ and ${\rm tr}(|D_k^{\kappa_1}D_\ell^{\kappa_2}\rho|) < \infty$ for all $\kappa_1+\kappa_2\le \kappa$, and thus
\[
C_{\rho,\kappa_1,\kappa_2} = {\rm tr}\big|[ X_2,\ldots[X_2,[X_1,\ldots[X_1,\rho]]]]\big| < \infty.
\]
But this implies that the density of ${\rm d}W_{h,\Phi}$ is contained in the Sobolev spaces $H^{p,\kappa}(\mathbb{R}^2)$ for all $2\le p\le \infty$.
\end{proof}

\begin{example}
Let $0< \lambda_1\le \lambda_2\le \cdots$ be an increasing sequence of positive numbers and $\{e_j|j\in\mathbb{N}\}$ a complete orthonormal system for $\eufrak{h}_\mathbb{C}$. Let $T_t:\eufrak{h}_\mathbb{C}\to\eufrak{h}_\mathbb{C}$ be the contraction semigroup defined by
\[
T_t e_j = e^{-t\lambda_j}e_j ,\qquad \text{ for } j\in\mathbb{N}, t\ge 0,
\]
with generator $A=\sum_{j\in\mathbb{N}} \lambda_j \mathbb{P}_j$. If the sequence increases fast enough to ensure that $\sum_{j=1}^\infty e^{-t\lambda_j}<\infty$, i.e.\ if ${\rm tr}\,T_t<\infty$ for $t>0$, then the second quantization $\rho_t=\Gamma(T_t):\eufrak{H}\to\eufrak{H}$ is a trace class operator with trace
\[
Z_t ={\rm tr}\, \rho_t = \sum_{\mathbf{n}\in\mathbb{N}_f^\infty} \langle e_{\mathbf{n}}, \rho_t e_{\mathbf{n}}\rangle
\]
where we use $\mathbb{N}_f^\infty$ to denote the finite sequences of non-negative integers and $\{e_\mathbf{n}|\mathbf{n}\in\mathbb{N}_f^\infty\}$ is the complete orthonormal system of $\eufrak{H}$ consisting of the vectors
\[
e_\mathbf{n}=e_1^{\odot n_1}\odot \cdots \odot e_r^{\odot n_r}, \qquad \text{ for } \mathbf{n}=(n_1,\ldots,n_r)\in \mathbb{N}_f^\infty,
\]
i.e.\ the symmetrization of the tensor $e_1\otimes \cdots \otimes e_1 \otimes \cdots\otimes e_r\otimes \cdots \otimes e_r$ where each vector $e_j$ appears $n_j$ times.
We get
$Z_t=\sum_{\mathbf{n}\in\mathbb{N}^\infty} \prod_{k=1}^\infty e^{-n_kt\lambda_k} = \prod_{k=1}^\infty \frac{1}{1-e^{-t\lambda_k}}$
for the trace of $\rho_t$. We shall be interested in the state defined by
\[
\Phi(X)=\frac{1}{Z_t}{\rm tr}(\rho_tX), \qquad\text{ for } X\in\mathcal{B}(\eufrak{H}).
\]
We get
\begin{eqnarray*}
\sum_{\mathbf{n}\in\mathbb{N}_f^\infty} \big|\langle e_\mathbf{n}, |\rho_{t/2} a^\ell(e_j)|^2e_\mathbf{n}\rangle\big| &=& \sum_{\mathbf{n}\in\mathbb{N}_f^\infty} || \rho_{t/2} a^\ell(e_j)||^2 \\
&=& \sum_{\mathbf{n}\in\mathbb{N}^\infty} n_j(n_j-1)\cdots (n_j-\ell+1) e^{-(n_j-\ell)t\lambda_j}\prod_{k\not=j} e^{-n_kt\lambda_k} \\
&\le & \sum_{n=0}^\infty (n+\ell)^{\ell} e^{-nt\lambda_j}   \prod_{k\not=j} \frac{1}{1-e^{-t\lambda_k}} < \infty,
\end{eqnarray*}
and therefore that $\rho_t a^\ell(e_j)$ defines a bounded operator with finite trace for all $j,\ell\in\mathbb{N}$ and $t>0$. Similarly we get
\[
{\rm tr}\, \big|a^\ell (e_j)\rho_t\big| < \infty, \quad
{\rm tr}\, \big|\rho_t \big(a^+(e_j)\big)^\ell\big| < \infty, \quad
\text{ etc. }
\]
and
\[
{\rm tr}\, \big|P^{\ell_1}(e_{j_1}) Q^{\ell_2}(e_{j_2})\rho_t\big| < \infty, \quad
{\rm tr}\, \big|P^{\ell_1}(e_{j_1}) \rho_t Q^{\ell_2}(e_{j_2})\big| < \infty, \quad \text{ etc. }
\]
for all $t>0$ and $j_1,j_2,\ell_1,\ell_2\in\mathbb{N}$. For a given $h\in\eufrak{h}\otimes\mathbb{R}^2$ with $\langle h_1,h_2\rangle\not=0$ (and thus in particular $h_1\not=0$ and $h_2\not=0$), we can always find indices $j_1$ and $j_2$ such that $\langle h_1,e_{j_1}\rangle \not=0$ and $\langle h_2,e_{j_2}\rangle \not=0$. Therefore it is not difficult to check that for all $\kappa\in\mathbb{N}$ all the conditions of Theorem \ref{theo regular2} are satisfied with $k=\left(\begin{array}{c}e_{j_1} \\ 0 \end{array}\right)$ and $\ell =\left(\begin{array}{c}0 \\ e_{j_2} \end{array}\right)$. We see that the Wigner density $w_{h,\Phi}$ of any pair $\big(P(h_1),Q(h_2)\big)$ with $\langle h_1,h_2\rangle\not=0$ in the state $\Phi(\cdot)= {\rm tr}(\rho_t\,\cdot)/Z_t$ is in $\bigcap_{\kappa\in\mathbb{N}} \bigcap_{2\le p\le\infty}H^{p,\kappa}(\mathbb{R}^2)$, in particular, its derivatives of all orders exist, and are bounded and square-integrable.
\end{example}

We will now show that this approach can also be applied to get sufficient conditions for the regularity of a single bounded self-adjoint operator, for simplicity we consider only vector states.

Given a bounded self-adjoint operator $X$, we call a measure $\mu_{X,\Phi}$ on the real line its distribution in the state $\Phi$, if all moments agree,
\[
\Phi(X^n) = \int_\mathbb{R} x^n{\rm d}\mu_{X,\Phi}, \qquad \text{ for all } n\in\mathbb{N}.
\]
Such a measure $\mu_{X,\Phi}$ always exists, it is unique and supported on the interval $\big[ -||X||,||X||\big]$.

\begin{proposition}\label{prop-reg-X}
Let $\omega\in\eufrak{H}$ be a unit vector and let $\Phi(\cdot)=\langle\omega,\cdot\,\omega\rangle$ be the corresponding vector state. The distribution $\mu_{X,\Phi}$ of an operator $X\in\mathcal{B}(\eufrak{H})$ in the state $\Phi$ has a bounded density, if there exists a $k\in\eufrak{h}_\mathbb{C}\otimes\mathbb{C}^2$ such that $\omega\in{\rm Dom}\big(Q(k_1)-P(k_2)\big)\cap{\rm Dom}\big(Q(\overline{k}_1)-P(\overline{k}_2)\big)$, $X\in{\rm Dom}\,D_k$, $X\cdot D_kX=D_kX\cdot X$, $D_kX$ invertible and $(D_kX)^{-1}\in {\rm Dom}\,D_k$.
\end{proposition}
\begin{proof}
Since $X\cdot D_kX=D_kX\cdot X$, we have $D_kp(X)= (D_kX)p'(X)$ for all polynomials $p$. We therefore get
\[
D_k\left( (D_kX)^{-1}p(X)\right) = p(X) D_k\left( (D_kX)^{-1}\right) + p'(X).
\]
The hypotheses of the proposition assure
\begin{eqnarray*}
\left|\langle\omega,D_k\left( (D_kX)^{-1}p(X)\right)\omega\rangle\right|&\!\le\!& \frac{||\big(\!Q(k_1)\!-\!P(k_2)\!\big)\omega||\!+\!||\big(\!Q(\overline{k}_1)\!-\!P(\overline{k}_2)\!\big)\omega||}{2}||(D_kX)^{-1}||\, ||p(X)|| \\
&\!\le\! &C_1 \sup_{x\in [ -||X||,||X||]} \big|p(x)\big|, \\[2mm]
\left|\langle\omega, p(X) D_k\left( (D_kX)^{-1}\right)\omega\rangle\right| 
&\!\le\! & \left|\left|D\left((D_kX)^{-1}\right)\right|\right|\, ||p(X)|| \\
&\!\le\! &C_2 \sup_{x\in [ -||X||,||X||]} \big|p(x)\big|,
\end{eqnarray*}
and therefore allow us to get the estimate
\begin{eqnarray*}
\left|\int_{-||X||}^{||X||} p'(x) {\rm d}\mu_{X,\Phi}(x)\right| &=&
\big|\langle\omega, p'(X)\omega\rangle\big| \\
&=& \Big|\Big\langle\omega,\left(D_k\Big( (D_kX)^{-1}p(X)\right)-p(X) D_k\left( (D_kX)^{-1}\right)\Big)\omega\Big\rangle\Big| \\
&\le & (C_1+C_2) \sup_{x\in [ -||X||,||X||]} \big|p(x)\big|
\end{eqnarray*}
for all polynomials $p$.
But this implies that $\mu_{X,\Phi}$ admits a bounded density.
\end{proof}

Let $n\in\mathbb{N}$, $n\ge 2$. We get that the density is even $n-1$ times differentiable, if in addition to the conditions of Proposition \ref{prop-reg-X} we also have
$X\in {\rm Dom}\,D_k^n$, $(D_kX)^{-1}\in {\rm Dom}\,D_k^n$, and
\[
\omega\in\bigcap_{1\le\kappa\le n}{\rm Dom}\big(Q(k_1)-P(k_2)\big)^\kappa \bigcap_{1\le\kappa\le n}{\rm Dom}\big(Q(\overline{k}_1)-P(\overline{k}_2)\big)^\kappa.
\]

The proof is similar to the proof of Proposition \ref{prop-reg-X}, we now use the formula
\[
p^{(n)}(X) = D_k^n\left( (D_kX)^{-n}p(X)\right) - \sum_{\kappa=0}^{n-1} A_\kappa p^{(\kappa)}(X),
\]
where $A_0,\ldots,A_{n-1}$ are some bounded operators, to get the necessary estimate
\[
\left|\int_{-||X||}^{||X||} p^{(n)}(x) {\rm d}\mu_{X,\Phi}(x)\right| \le
C \sup_{x\in [ -||X||,||X||]} \big|p(x)\big|
\]
by induction over $n$.

\subsection{The white noise case}\label{white noise}

Let now $\eufrak{h}=L^2(T,\mathcal{B},\mu)$, where $(T,\mathcal{B},\mu)$ is a measure space such that $\mathcal{B}$ is countably generated. In this case we can apply the divergence operator to processes indexed by $T$, i.e.\ $\mathcal{B}(\eufrak{H})$-valued measurable functions on $T$, since they can be interpreted as elements of the Hilbert module, if they are square-integrable. Let $L^2\big(T,\mathcal{B}(\eufrak{H})\big)$ denote all $\mathcal{B}(\eufrak{H})$-valued measurable functions $t\mapsto X_t$ on $T$ with $\int_T ||X_t||^2 {\rm d}t< \infty$.
Then the definition of the divergence operator becomes
\begin{eqnarray*}
{\rm Dom}\,\delta &=& \Big\{ X=(X^1,X^2)\in L^2\big(T,\mathcal{B}(\eufrak{H})\big)\oplus L^2\big(T,\mathcal{B}(\eufrak{H})\big)\Big| \\
&&\exists M\in\mathcal{B}(\eufrak{H}), \forall k_1,k_2\in \eufrak{h}_\mathbb{C}: \langle\mathcal{E}(k_1), M\mathcal{E}(k_2)\rangle \\
&&= \int_T\big(i(k_2-\overline{k}_1)\langle\mathcal{E}(k_1), X^1(t)\mathcal{E}(k_2)\rangle+(\overline{k}_1+k_2) \langle\mathcal{E}(k_1), X^2(t)\mathcal{E}(k_2)\rangle\big){\rm d}\mu(t)\Big\},
\end{eqnarray*}
and $\delta(X)$ is  equal to the unique operator satisfying the above condition. We will also use the notation
\[
\delta(X)=\int_T X^1(t){\rm d}P(t) + \int_T X^2(t){\rm d}Q(t),
\]
and call $\delta(X)$ the Hitsuda-Skorohod integral of $X$.

Belavkin \cite{belavkin91a,belavkin91b} and Lindsay \cite{lindsay93} have defined non-causal quantum stochastic integrals with respect to the creation, annihilation, and conservation processes on the boson Fock space over $L^2(\mathbb{R}_+)$ using the classical derivation and divergence operators. It turns out that our Hitsuda-Skorohod integral coincides with their creation and annihilation integrals, up to a coordinate transformation. This immediately implies that for adapted, integrable processes our integral coincides with the quantum stochastic creation and annihilation integrals defined by Hudson and Parthasarathy, cf.\ \cite{hudson+parthasarathy84,parthasarathy92}.

Let us briefly recall, how they define the creation and annihilation integral, cf.\ \cite{lindsay93}. They use the derivation operator $\tilde{D}$ and the divergence operator $\tilde{\delta}$ from the classical calculus on the Wiener space $L^2(\Omega)$, defined on the Fock space $\Gamma(L^2(\mathbb{R}_+;\mathbb{C}))$ over $L^2(\mathbb{R}_+;\mathbb{C})=L^2(\mathbb{R}_+)_\mathbb{C}$ via the isomorphism between $L^2(\Omega)$ and $\Gamma\big(L^2(\mathbb{R}_+;\mathbb{C})\big)$. On the exponential vectors $\tilde{D}$ acts as
\[
\tilde{D}\mathcal{E}(k)=\mathcal{E}(k)\otimes k, \qquad \text{ for } k\in L^2(\mathbb{R}_+,\mathbb{C}),
\]
and $\tilde{\delta}$ is its adjoint. Note that due to the isomorphism between $\Gamma\big(L^2(\mathbb{R}_+;\mathbb{C})\big)\otimes L^2(\mathbb{R}_+;\mathbb{C})$ and $L^2\Big(\mathbb{R}_+;\Gamma\big(L^2(\mathbb{R}_+;\mathbb{C})\big)\Big)$, the elements of  $\Gamma\big(L^2(\mathbb{R}_+;\mathbb{C})\big)\otimes L^2(\mathbb{R}_+;\mathbb{C})$ can be interpreted as function on $\mathbb{R}_+$. In particular, for the exponential vectors we get $\big(D\mathcal{E}(k)\big)_t=k(t)\mathcal{E}(k)$ almost surely. The action of the annihilation integral $\int F_t dA_t$ on some vector $\psi\in\Gamma\big(L^2(\mathbb{R}_+;\mathbb{C})\big) $ is then defined as the Bochner integral
\[
\int^{\rm BL}_{\mathbb{R}_+} F_t dA_r\psi = \int_{\mathbb{R}_+} F_t (D\psi)_t {\rm d}t,
\]
and that of the creation integral as
\[
\int^{\rm BL}_{\mathbb{R}_+} F_t {\rm d}A^*_t\psi = \tilde{\delta}(F_\cdot\psi).
\]
These definitions satisfy the adjoint relations
\[
\left(\int^{\rm BL}_{\mathbb{R}_+} F_t dA_r\right)^* \supset \int^{\rm BL}_{\mathbb{R}_+} F^*_t {\rm d}A^*_t, \quad \text{ and } \quad \left(\int^{\rm BL}_{\mathbb{R}_+} F_t dA^*_r\right)^* \supset \int^{\rm BL}_{\mathbb{R}_+} F^*_t {\rm d}A_t.
\]

\begin{proposition}
Let $(T,\mathcal{B},\mu)=(\mathbb{R}_+,\mathcal{B}(\mathbb{R}_+),{\rm d}x)$, i.e.\ the positive half-line with the Lebesgue measure, and let $X=(X^1,X^2)\in {\rm Dom}\,\delta$. Then we have
\[
\int_{\mathbb{R}_+} X^1(t){\rm d}P(t) + \int_{\mathbb{R}_+} X^2(t){\rm d}Q(t) = \int^{\rm BL}_{\mathbb{R}_+}(X^2-iX^1) {\rm d}A^*_t + \int^{\rm BL}_{\mathbb{R}_+}(X^2+iX^1){\rm d}A_t.
\]
\end{proposition}
\begin{proof}
To prove this, we show that the Belavkin-Lindsay integrals satisfy the same formula for the matrix elements between exponential vectors. Let $(F_t)_{t\in\mathbb{R}_+}\in L^2\big(\mathbb{R}_+,\mathcal{B}(\eufrak{H})\big)$ be such that its creation integral in the sense of Belavkin and Lindsay is defined with a domain containing the exponential vectors. Then we get
\begin{eqnarray*}
\langle\mathcal{E}(k_1),\int^{\rm BL}_{\mathbb{R}_+} F_t {\rm d}A^*_t\mathcal{E}(k_2)\rangle &=& \langle\mathcal{E}(k_1), \tilde{\delta}\big(F_\cdot\mathcal{E}(k_2)\big)\rangle \\
&=& \langle\big(\tilde{D}\mathcal{E}(k_1)\big)_\cdot, F_\cdot\mathcal{E}(k_2)\rangle \\
&=& \int\overline{k_1(t)} \langle\mathcal{E}(k_1), F_t\mathcal{E}(k_2)\rangle{\rm d}t.
\end{eqnarray*}
For the annihilation integral we deduce the formula
\begin{eqnarray*}
\langle\mathcal{E}(k_1),\int^{\rm BL}_{\mathbb{R}_+} F_t {\rm d}A_t\mathcal{E}(k_2)\rangle &=& \langle\int^{\rm BL}_{\mathbb{R}_+} F^*_t {\rm d}A^*_t\mathcal{E}(k_1),\mathcal{E}(k_2)\rangle \\
&=& \overline{\int_{\mathbb{R}_+}\overline{k_2(t)} \langle\mathcal{E}(k_2), F^*_t\mathcal{E}(k_1)\rangle{\rm d}t} \\
&=&\int_{\mathbb{R}_+} k_2(t) \langle\mathcal{E}(k_1), F_t\mathcal{E}(k_2)\rangle{\rm d}t.
\end{eqnarray*}
\end{proof}

The integrals defined by Belavkin and Lindsay are an extension of those defined by Hudson and Parthasarathy, therefore the following is now obvious.

\begin{corollary}
For adapted processes $X\in{\rm Dom}\,\delta$, the Hitsuda-Skorohod integral
\[
\delta(X)=\int_T X^1(t){\rm d}P(t) + \int_T X^2(t){\rm d}Q(t)
\]
coincides with the Hudson-Parthasarathy quantum stochastic integral defined in \cite{hudson+parthasarathy84}.
\end{corollary}

\subsection{Iterated integrals}

We give a short, informal discussion of iterated integrals of deterministic functions, showing a close relation between these iterated integrals and the so-called Wick product or normal-ordered product. Doing so, we will encounter unbounded operators, so that strictly speaking we have not defined $\delta$ for them. But everything could be made rigorous by choosing an appropriate common invariant domain for these operators, e.g., vectors with a finite chaos decomposition.

In order to be able to iterate the divergence operator, we define $\delta$ on $\mathcal{B}(\eufrak{H})\otimes \eufrak{h}_\mathbb{C}\otimes\mathbb{C}^2\otimes H$, where $H$ is some Hilbert space, as $\delta\otimes {\rm id}_{H}$.

Using Equation (\ref{wick}), on can show by induction
\[
\delta^n \left(\begin{array}{c} h^{(1)}_1 \\ h^{(1)}_2 \end{array}\right) \otimes \cdots \otimes \left(\begin{array}{c} h^{(n)}_1 \\ h^{(n)}_2 \end{array}\right)
= \sum_{I\subseteq\{1,\ldots,n\}} \prod_{j\in I} a^+(h^{(j)}_2-i h^{(j)}_1) \prod_{j\in \{1,\ldots,n\}\backslash I} a(\overline{h}^{(j)}_2+i \overline{h}^{(j)}).
\]
for $h^{(1)}=\left(\begin{array}{c} h^{(1)}_1 \\ h^{(1)}_2 \end{array}\right),\ldots,h^{(n)}=\left(\begin{array}{c} h^{(n)}_1 \\ h^{(n)}_2 \end{array}\right)\in\eufrak{h}_\mathbb{C}\otimes\mathbb{C}^2$. This is just the Wick product of $\big(P(h^{(1)}_1+Q(h^{(1)}_2)\big),\ldots,\big(P(h^{(n)}_1+Q(h^{(n)}_2)\big)$, i.e.
\[
\delta^n \left(\begin{array}{c} h^{(1)}_1 \\ h^{(1)}_2 \end{array}\right) \otimes \cdots \otimes \left(\begin{array}{c} h^{(n)}_1 \\ h^{(n)}_2 \end{array}\right)
= \big(P(h^{(1)}_1+Q(h^{(1)}_2)\big)\diamond\cdots\diamond\big(P(h^{(n)}_1+Q(h^{(n)}_2)\big),
\]
where the Wick product $\diamond$ is defined on the algebra generated by $\{P(k),Q(k)|k\in\eufrak{h}_\mathbb{C}\}$ by
\begin{eqnarray*}
P(h)\diamond X &=& X\diamond P(h) = -ia^+(h)X +iX a(\overline{h}) \\
Q(h)\diamond X &=& X\diamond Q(h) = a^+(h)X+ Xa(\overline{h})
\end{eqnarray*}
for $X\in {\rm alg}\,\{P(k),Q(k)|k\in\eufrak{h}_\mathbb{C}\}$ and $h\in\eufrak{h}_\mathbb{C}$ in terms of the momentum and position operators, or, equivalently, by
\begin{eqnarray*}
a^+(h)\diamond X &=& X\diamond a^+(h) =  a^+(h)X \\
a(h)\diamond X &=& X\diamond a(h) = Xa(h)
\end{eqnarray*}
in terms of creation and annihilation.

\section{Conclusion}

We have defined a derivation operator $D$ and a divergence operator $\delta$ on $\mathcal{B}(\eufrak{H})$ and $\mathcal{B}(\eufrak{H},\eufrak{H}\otimes\eufrak{h}_\mathbb{C}\otimes\mathbb{C}^2)$, resp., and shown that they have similar properties as the derivation operator and the divergence operator in classical Malliavin calculus. As far as we know, this is the first time that $D$ and $\delta$ are considered as operators defined on a non-commutative operator algebra, except for the free case \cite{biane+speicher99}, where the operator algebra is isomorphic to the full Fock space. To obtain close analogues of the classical relations involving the divergence operator, we needed to impose additional commutativity conditions, but Proposition \ref{prop no go} shows that this can not be avoided. Also, its domain is rather small, because we require $\delta(u)$ to be a bounded operator and so, e.g., the ``deterministic'' elements $h\in\eufrak{h}_\mathbb{C}\otimes\mathbb{C}^2$ are not integrable unless $h=0$. One of the main goals of our approach is the study of Wigner densities as joint densities of non-commutating random variables. We showed that the derivation operator can be used to obtain sufficient conditions for its regularity, see, e.g., Theorem \ref{theo regular2}. It seems likely that these results can be generalized by weakening or modifying the hypotheses. It would be interesting to apply these methods to quantum stochastic differential equations.

\section*{Acknowledgements}

One of the authors (U.F.) wants to thank K.B.~Sinha and B.V.R.~Bhat of the Indian Statistical Institutes in Bangalore and Dehli, where parts of this work were carried out, for their generous hospitality, and the DFG for a travel grant. R.S.\ and U.F.\ are also grateful to the Institut Henri Poincar\'e in Paris for kind hospitality during the special semester on free probability and operator spaces. U.F.\ is also thankful to Michael Skeide for introducing him to the basics of Hilbert module theory.


\end{document}